\documentclass{amsart}
\usepackage{amsmath}
  \usepackage{paralist}
  \usepackage{graphics} 
  \usepackage{epsfig} 
\usepackage{graphicx}  \usepackage{epstopdf}
 \usepackage[colorlinks=true]{hyperref}
 \usepackage{amsmath}
\usepackage{amsfonts}
\usepackage{amssymb}
\usepackage{ mathrsfs }
\usepackage{amsthm}
\usepackage[pagewise]{lineno}
\usepackage{enumitem}
\usepackage{graphicx,color,xcolor,tikz}
\usepackage{bm}
\usepackage{hyperref}
\hypersetup{
    colorlinks=true,                         
    linkcolor=blue, 
    citecolor=red, 
  } 
\hypersetup{urlcolor=blue, citecolor=red}

\newtheorem{thm}{Theorem}
\newtheorem{prop}[thm]{Proposition}

\newtheorem{defi}[thm]{Definition}
\newtheorem{rem}[thm]{Remark}

\def\A{\mathcal A}
\def\B{\mathcal B}

\def\I{\mathcal I}
\def\M{\mathcal M}

\def\R{\mathbb R}

\def\W{\mathcal W}
\def\e{\varepsilon}

\newcommand{\pt}{\partial}

\newcommand{\abs}[1]{\ensuremath{\left|#1\right|}}

\newcommand{\norm}[1]{\ensuremath{\left\|#1\right\|}}
\makeatletter
\newcommand{\doublewidetilde}[1]{{%
  \mathpalette\double@widetilde{#1}%
}}
\newcommand{\double@widetilde}[2]{%
  \sbox\z@{$\m@th#1\widetilde{#2}$}%
  \ht\z@=.9\ht\z@
  \widetilde{\box\z@}%
}
\makeatother

\title[Derivation of a new macroscopic bidomain model including three scales]{Derivation of a new macroscopic bidomain model including three scales for the electrical activity of cardiac tissue}


\author{Fakhrielddine Bader$^*$ }
\address[Fakhrielddine Bader]{Mathematics Laboratory, Doctoral school of Sciences and Technologies, Lebanese University, Hadat Beirut, Lebanon \&  Laboratoire de Mathématiques Jean Leray, École Centrale de Nantes, Nantes, France}
\email{fakhrielddine.bader@ec-nantes.fr}

\author{Mostafa Bendahmane}
\address[Mostafa Bendahmane]{Institut de Mathématiques de Bordeaux and INRIA-Carmen Bordeaux Sud-Ouest, Université de Bordeaux, 33076 Bordeaux Cedex, France}
\email{mostafa.bendahmane@u-bordeaux.fr}

\author{Mazen Saad}
\address[Mazen Saad]{Laboratoire de Mathématiques Jean Leray, UMR 6629 CNRS, École Centrale de Nantes,  1 rue de Noé, 44321 Nantes, France}
\email{mazen.saad@ec-nantes.fr}

\author{Raafat Talhouk}
\address[Raafat Talhouk]{Mathematics Laboratory, Doctoral school of Sciences and Technologies, Lebanese University, Hadat Beirut, Lebanon}
\email{rtalhouk@ul.edu.lb}

\thanks{$^*$ Corresponding author: fakhri.bader.fb@gmail.com}

\begin{document}
\maketitle



\begin{abstract}
In the present paper, a new three-scale asymptotic homogenization method is proposed to study the electrical behavior of the cardiac tissue structure with multiple heterogeneities at two different levels. The first level is associated with the mesoscopic structure such that the cardiac tissue is composed of extracellular and intracellular domains. The second level is associated with the microscopic structure in such a way the intracellular medium can only be viewed as a periodical layout of unit cells (mitochondria). Then, we define two kinds of local cells that are obtained by upscaling methods. The homogenization method is based on a power series expansion which allows determining the macroscopic (homogenized) bidomain model from the microscopic bidomain problem at each structural level. First, we use the two-scale asymptotic expansion to homogenize the extracellular problem. Then, we apply a three-scale asymptotic expansion in the intracellular problem to obtain its homogenized equation at two levels. The first upscaling level of the intracellular structure yields the mesoscopic equation and the second step of the homogenization leads to obtain the intracellular homogenized equation. Both the mesoscopic and microscopic information is obtained by homogenization to capture local characteristics inside the cardiac tissue structure. Finally, we obtain the macroscopic bidomain model and the heart domain coincides with the intracellular medium and extracellular one, which are two inter-penetrating and superimposed continua connected at each point by the cardiac cellular membrane. The interest of the proposed method comes from the fact that it combines microscopic and mesoscopic characteristics to obtain a macroscopic description of the electrical behavior of the heart.
\end{abstract}
\tableofcontents

\section{Introduction}
The heart is an organ that ensures life for all living beings. Indeed, its great importance comes from its organic function which allows the circulation of blood throughout the body. It is a muscular organ composed of four cavities: the left atrium and ventricle which represent the left heart, the right atrium and ventricle which form the right heart. These four cavities are surrounded by a cardiac tissue that is organized into muscle fibers. These fibers form a network of cardiac muscle cells called "cardiomyocytes" connected end-to-end by junctions. For more details about the physiological background, we refer to \cite{Keener} and about the electrical activity of the heart we refer to \cite{char}.
 
\begin{figure}[h!]
  \centering
  \includegraphics[width=13cm]{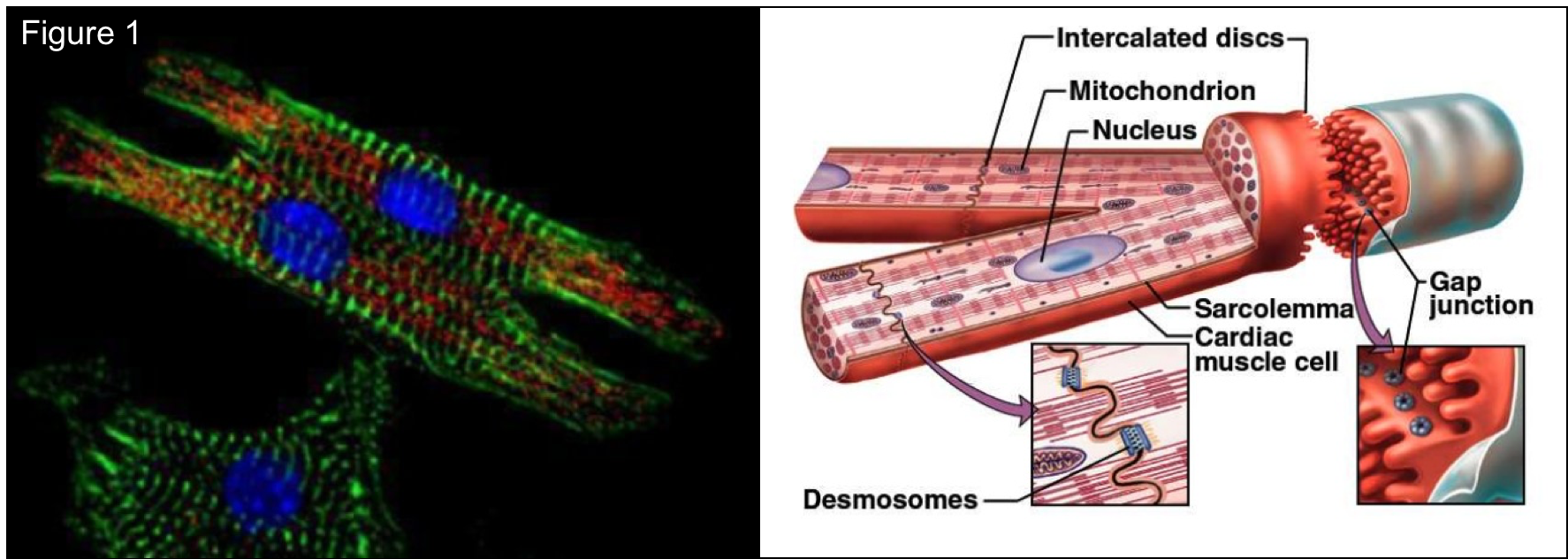}
  \caption{Representation of the cardiomyocyte structure}
 \url{http://www.cardio-research.com/cardiomyocytes}
  \label{cardio}
 \end{figure}

 The structure of cardiac tissue (myocarde) studied in this paper is characterized at three different scales (see Figure \ref{cardio}). At mesoscopic scale, the cardiac tissue is divided into two media: one contains the contents of the cardiomyocytes, in particular the "cytoplasm" which is called the "intracellular" medium, and the other is called extracellular and consists of the fluid outside the cardiomyocytes cells. These two media are separated by a cellular membrane (the sarcolemma) allowing the penetration of proteins, some of which play a passive role and others play an active role powered by cellular metabolism. At microscopic scale, the cytoplasm comprises several organelles such as mitochondria. Mitochondria are often described as the "energy powerhouses" of cardiomyocytes and are surrounded by another membrane. Then, we consider only that the intracellular medium can be viewed as a periodic perforated structure composed of other connected cells. While at the macroscopic scale, this domain coincides with the intracellular medium  and extracellular  one, which are two inter-penetrating and superimposed continua connected at each point by the cardiac cellular membrane.
 
 It should be noted that there is a difference between the chemical composition of the cytoplasm and that of the extracellular medium. This difference plays a very important role in cardiac activity. In particular, the concentration of anions (negative ions) in cardiomyocytes is higher than in the external environment. This difference of concentrations creates a transmembrane potential, which is the difference in potential between these two media. The  model that describes the electrical activity of the heart, is called by "Bidomain model". The first mathematical formulation of this model was constructed by Tung \cite{tung1978}. The authors in \cite{bendunf19} established the well-posedness of this microscopic bidomain model under different conditions and proved the existence and uniqueness of their solutions (see other work \cite{marco}). \\
The microscopic bidomain model \cite{henri,colli05} consists of two quasi-static equations, one for the electrical potential in the intracellular medium and one for the extracellular medium, coupled through a dynamic boundary equation at the interface of the two regions (the membrane $\Gamma^{y}$). In our study, these equations depend on two small scaling parameters $\e$ and $\delta$ whose are respectively the ratio of the microscopic and mesoscopic scales from the macroscopic scale. 
  
   Our goal in this paper is to derive, using a homogenization method, the macroscopic (homogenized) bidomain model of the cardiac tissue from the microscopic bidomain model. In general, the homogenization theory is the analysis of the macroscopic behavior of biological tissues by taking into account their complex microscopic structure. For an introduction to this theory, we cite \cite{sanchez}, \cite{doina},\cite{tartarintro} and \cite{bakhvalov}. Applications of this technique can be found in modeling solids, fluids, solid-fluid interaction, porous media, composite materials, cells, and cancer invasion. 
Several methods are related to this theory. First, the multiple-scale method established by Benssousan et al. \cite{ben} and used in mechanics and physics for problems containing several small scaling parameters.
   The two-scale convergence method was introduced by Nugesteng \cite{ngu} and developed by Allaire et al. \cite{allaire92}. In addition, Allaire et al. \cite{allairebriane}, Trucu et al. \cite{trucu} introduced a further generalization of the previous method via a three-scale convergence approach for distinct problems. Here we are not dealing with a rigorous multi-scale convergence setting, as our main motivation lies in the direct application of asymptotic homogenization by following a formal approach and accounting for a novel series expansion in terms of two distinct scaling parameters $\e$ and $\delta.$ Recently, the periodic unfolding method introduced by Cioranescu et al. for the fixed domains in \cite{doinaunf02}  and for the perforated domains in \cite{doinaunf06}.  This method is essentially based on two operators: the first represents the unfolding operator and the second operator consists to separate the microscopic and macroscopic scales 
   (see also \cite{doinaunf08,doinaunf12}).

  There are some of these methods that are applied on the microscopic bidomain model to obtain the homogenized macroscopic model. First, Krassowska and Neu \cite{neukra} applied the two-scale asymptotic method to formally obtain this macroscopic model (see also \cite{henri,richardson} for different approaches). Furthermore, Pennachio et al. \cite{colli05} used the tools of the $\Gamma$-convergence method to obtain a rigorous mathematical form of this homogenized macroscopic model.  In \cite{anna,karlsen}, the authors used the theory of two-scale convergence method to derive the homogenized bidomain model. Recently, the authors in \cite{bendunf19} proved the existence and uniqueness of solution of the microscopic bidomain model by using the Faedo-Galerkin method. Further, they applied the unfolding homogenization method at two scales. Some recent works are available on the numerical implementation of bidomain models in the context of pure electro-physiology in \cite{kaliske2012} and of cardiac electromechanics in \cite{kaliske2013,kaliske2018}.
  
  
  \textit{The main of contribution of the present paper.} The cardiac tissue structure viewed at micro-macro scales and studied at the three different scales where the intracellular medium is a periodic composed of connected cells. The aim is to derive the two-levels homogenized bidomain model of cardiac electro-physiology from the microscopic bidomain model. This paper presents a formal mathematical writing for the results obtained in a recent work \cite{BaderUnf} based on a three-scale unfolding homogenization method. In \cite{BaderUnf}, we used unfolding operators, to converge our meso-microscopic bidomain problem as $\e,\delta\rightarrow 0$ and then to obtain the same macroscopic bidomain system. While in the present work, we will apply the two-scale asymptotic expansion method on the extracellular medium (similar derivation may be found in \cite{henri}). Further, we will derive a formal approach, by accounting for a three-scale asymptotic expansion, in terms of two distinct scaling parameters $\e$ and $\delta$ on the intracellular medium based on the work of Benssousan et al. \cite{ben}. The asymptotic expansion proposed to investigate the effective properties of the cardiac tissue at each structural level, namely, micro-meso-macro scales. Moreover, to treat the bidomain problem in this work, the multi-scale technique is needed to be established in time domain directly.  
   
  
  \textit{The outline of the paper is as follows.} In Section \ref{geobid}, we introduce the microscopic bidomain model in the cardiac tissue structure featured by two parameters, $\e$ and $\delta,$ characterizing the  micro and mesoscopic scales. Section \ref{asym} is devoted to homogenization procedure. The two-scale asymptotic expansion method applied in the extracellular problem is explained in Subsection \ref{asymextra}. The homogenized equation for the extracellular problem are obtained in terms of the coefficients of conductivity matrices and cell problems. In Subsection \ref{asymintra}, the homogenized equation for the intracellular problem is obtained similarly at two levels but using a three-scale asymptotic expansion which depends on $\e$ and $\delta.$ The first level of homogenization yields the mesoscopic problem and then we complete the second level to obtain the corresponding homogenized equation. Finally, the main result is presented in Section \ref{macro_bid_asym} and the macroscopic bidomain equations are recuperated from the extracellular and intracellular homogenized equations. In Appendix \ref{appA}, we report some notations and special functional spaces used for the homogenization. More properties and theorems including these spaces are also postponed in the Appendix \ref{appA}.

\section{Bidomain modeling of the heart tissue}\label{geobid}
The aim of this section is to define the geometry of cardiac tissue and to present the microscopic bidomain model of the heart.
\subsection{Geometric idealization of the myocardium micro-structure}
\begin{figure}[h!]
  \centering
  \includegraphics[width=15cm]{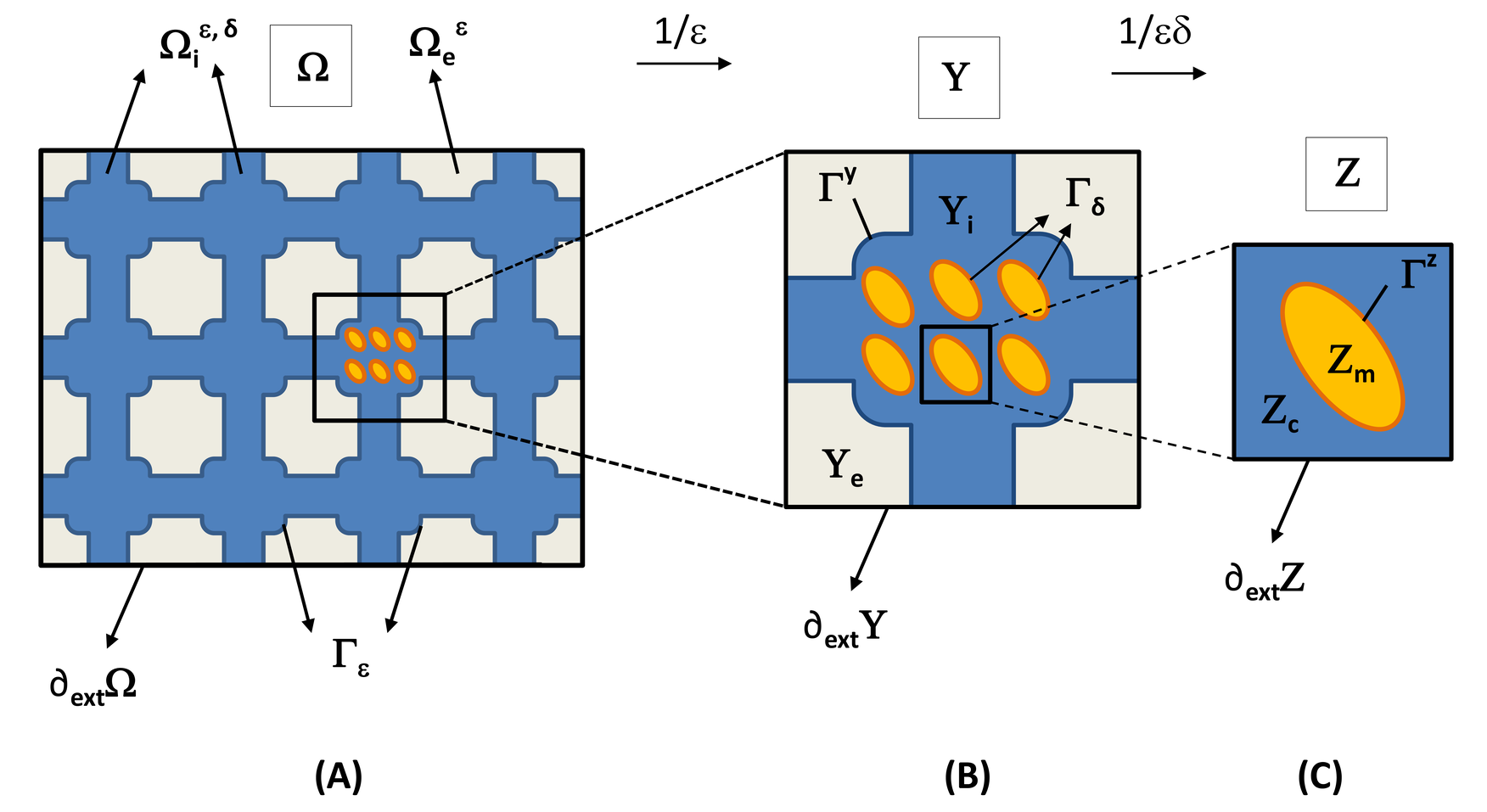}
  \caption{(A) Periodic heterogeneous domain $\Omega$, (B) Unit cell $Y$ at $\e$-structural level and (C) Unit cell $Z$ at $\delta$-structural level }
  \label{three_scale}
 \end{figure}

The cardiac tissue $\Omega \subset \R^d$ is considered as a heterogeneous double periodic domain with a Lipschitz boundary $\pt \Omega$. The structure of the tissue is periodic at mesoscopic and microscopic scales related to two small parameters $\e$ and $\delta$, respectively, see Figure \ref{three_scale}. 

 Following the standard approach of the homogenization theory, this structure is featured by two parameters $\ell^\text{mes}$ and $\ell^\text{mic}$ characterizing, respectively, the mesoscopic and microscopic length of a cell. Under the two-level scaling, the characteristic lengths $\ell^\text{mes}$ and $\ell^\text{mic}$ are related to a given macroscopic  length $L$ (of the cardiac fibers), such that the two scaling parameters $\e$ and $\delta$ are introduced by:
 $$\e=\frac{\ell^\text{mes}}{L} \text{ and } \delta=\frac{\ell^\text{mic}}{L} \ \text{ with } \ \ell^\text{mic}<<\ell^\text{mes}.$$

\paragraph{\textbf{The mesoscopic scale.}} The domain $\Omega$ is composed of two ohmic volumes, called intracellular $\Omega_i^{\e,\delta}$ and extracellular $\Omega_e^{\e}$ medium (for more details see \cite{colli05}). Geometrically, we find that $\Omega_i^{\e,\delta}$ and $\Omega_e^{\e}$ are two open connected regions such that:
 $$\overline{\Omega}=\overline{\Omega}_{i}^{\e,\delta}\cup \overline{\Omega}_e^{\e}, \ \text{with} \ \Omega_i^{\e,\delta}\cap\Omega_e^{\e}=\emptyset.$$ These two regions are separated by the surface membrane $\Gamma_{\e}$ which is expressed by: $$\Gamma_{\e}=\pt \Omega_i^{\e,\delta} \cap \pt \Omega_e^{\e},$$ assuming that the membrane is regular. We can observe that the domain $\Omega_i^{\e,\delta}$ as a perforated domain obtained from $\Omega$ by removing the holes which correspond to the extracellular domain $\Omega_e^{\e}.$ 
 
  At this $\bm{\varepsilon}$-structural level, we can divide $\Omega$ into $N_{\e}$ small elementary cells $Y_{\e}=\overset{d}{\underset{n=1}{\prod }}]0,\e \, \ell^\text{mes}_n[,$  with $\ell^\text{mes}_1,\dots,\ell^\text{mes}_d$ are positive numbers. These small cells are all equal, thanks to a translation and scaling by $\e,$ to the same unit cell of periodicity called the reference cell  $Y=\overset{d}{\underset{n=1}{\prod }}]0,\ell^\text{mes}_n[.$ Indeed, if we denote by $T_\e^k$ a translation of $\e k$   with $k=( k_1,\dots, k_d ) \in\mathbb{Z}^d.$ In addition, if the cell considered $Y^k_\e$ is located at the $k_i^{\text{\text{ième}}}$ position according to the direction $i$ of space considered, we can write:
\begin{equation*}
 Y^k_{\e}:=T^k_\e+\e Y=\lbrace \e \xi : \xi \in k_\ell+Y \rbrace,
 \end{equation*}
 with $ k_\ell:=( k_1\ell^\text{mes}_1,\dots,  k_d \ell^\text{mes}_d ).$\\
 Therefore, for each macroscopic variable $x$ that belongs to $\Omega,$ we define the corresponding mesoscopic variable $y\approx\dfrac{x}{\e}$ that belongs to $Y$ with a  translation. Indeed, we have:
 \begin{equation*}
 x \in \Omega \Rightarrow \exists k \in\mathbb{Z}^d  \ \text{ such that }  \ x \in Y^k_\e \Rightarrow x=\e (k_\ell+y) \Rightarrow y=\dfrac{x}{\e}-k_\ell \in Y.
 \end{equation*}

 Since, we will study in the extracellular medium $\Omega^{\e}_e$ the behavior of the functions $u(x,y)$ which are $\textbf{y}$-periodic, so by periodicity we have $u\left( x,\dfrac{x}{\e}-k_\ell\right)  =u\left( x,\dfrac{x}{\e}\right) .$ By notation, we say that $y=\dfrac{x}{\e}$ belongs to $Y.$

 We are assuming that the cells are periodically organized as a regular network of interconnected cylinders at the mesoscale. The mesoscopic unit cell $Y$ is also divided into two parts: intracellular $Y_i$ and extracellular $Y_e.$ These two parts are separated by a common boundary $\Gamma^{y}.$ So, we have:
 \begin{equation*}
 Y=Y_i \cup Y_e \cup \Gamma^{y}, \quad \Gamma^{y}= \pt Y_i \cap \pt Y_e.
 \end{equation*}  
In a similar way, we can write the corresponding common periodic boundary as follows:
  \begin{equation*}
 \Gamma^k_{\e}:=T^k_\e+\e \Gamma^{y}=\lbrace \e \xi : \xi \in k_\ell+\Gamma^{y} \rbrace,
 \end{equation*}
 with   $T^k_\e$ denote the same previous translation.

 In summary, the intracellular and extracellular medium can be described as the intersection of the cardiac tissue $\Omega$ with the cell $Y^k_{j,\e}$ for $j=i,e:$  
 \begin{equation*}
\Omega^{\e,\delta}_i=\Omega \cap \underset{k\in \mathbb{Z}^d}{\bigcup} Y^k_{i,\e}, \quad \Omega^{\e}_e=\Omega \cap \underset{k\in \mathbb{Z}^d}{\bigcup} Y^k_{e,\e}, \quad \Gamma_{\e}=\Omega \cap \underset{k\in \mathbb{Z}^d}{\bigcup} \Gamma^k_{\e},
\end{equation*}
with each cell defined by $Y^k_{j,\e}=T^k_\e+\e Y_j $ for $j=i,e$.

\paragraph{\textbf{The microscopic scale.}} The cytoplasm  contain far more mitochondria described as "the powerhouse of the myocardium" surrounded by another membrane $\Gamma_{\delta}.$ then, we only assume that the intracellular medium $\Omega_i^{\e,\delta}$ can also be viewed as a periodic perforated domain.
 
 At this $\bm{\delta}$-structural level, we can divide this medium with the same strategy into small elementary cells $Z_{\delta}=\overset{d}{\underset{n=1}{\prod }}]0,\delta \, \ell^\text{mic}_n[,$ with $\ell^\text{mic}_1,\dots,\ell^\text{mic}_d$ are positive numbers. Using a similar translation (noted by $T^{k'}_\delta$), we return to the same reference cell noted by $Z=\overset{d}{\underset{n=1}{\prod }}]0,\ell^\text{mic}_n[.$ Note that if the cell considered $Z^{k'}_{\delta}$ is located at the $k_n^{'\text{ième}}$ position according to the direction $n$ of space considered, we can write:
\begin{equation*}
 Z^{k'}_{\delta}:=T^{k'}_\delta+\delta Z=\lbrace \delta \zeta : \zeta \in k'_{\ell'}+Z \rbrace,
 \end{equation*}
 with $ k'_{\ell'} :=( k'_1\ell^\text{mic}_1,\dots,  k'_d \ell^\text{mic}_d ).$ \\
 Therefore, for each macroscopic variable $x$ that belongs to $\Omega,$ we also define the corresponding microscopic variable $z\approx\dfrac{y}{\delta}\approx\dfrac{x}{\e\delta}$ that belongs to $Z$ with the translation $T^{k'}_\delta$.
 
The microscopic reference cell $Z$ splits into two parts: mitochondria part $Z_{m}$ and the complementary part $Z_{c}:=Z \setminus Z_{m}.$ These two parts are separated by a common boundary $\Gamma^{z}.$ So, we have:
 \begin{equation*}
 Z=Z_{m} \cup Z_{c} \cup \Gamma^{z}, \quad \Gamma^{z}= \pt Z_{m}.
 \end{equation*}
By definition, we have $\pt Z_{c}=\pt_{\text{ext}} Z \cup \Gamma^{z}.$
 
  More precisely, we can write the intracellular meso- and microscopic domain $\Omega_i^{\e,\delta}$ as follows:
\begin{equation*}
 \Omega^{\e,\delta}_i=\Omega \cap \underset{k\in \mathbb{Z}^d}{\bigcup} \left( Y^k_{i,\e} \cap \underset{k'\in \mathbb{Z}^d}{\bigcup} Z^{k'}_{c,\delta}\right) 
 \end{equation*}
with $Z^{k'}_{c,\delta}$ defined by:
\begin{equation*}
 Z^{k'}_{c,\delta} :=T^{k'}_\delta+\delta Z_{c}=\lbrace \delta \zeta : \zeta \in k'_{\ell'}+Z_{c} \rbrace.
 \end{equation*}
  In the intracellular medium $\Omega^{\e,\delta}_i,$ we will study the behavior of the functions  $u(x,y,z)$ which are $\textbf{z}$-periodic, so by periodicity we have $u\left( x,y,\dfrac{x}{\e \delta}-\dfrac{k_{\ell}}{\delta}-k'_{\ell'}\right) =u\left( x,y,\dfrac{x}{\e\delta}\right).$ By notation, we say that $z=\dfrac{x}{\e\delta}$ belongs to $Z.$ Similarly, we describe the common boundary at microscale as follows: 
 \begin{equation*}
\Gamma_{\delta}=\Omega \cap \underset{k'\in \mathbb{Z}^d}{\bigcup} \Gamma^{k'}_{\delta},
\end{equation*}
where $\Gamma^{k'}_{\delta}$ is given by:
\begin{equation*}
 \Gamma^{k'}_{\delta}:=T^{k'}_\delta+\delta \Gamma^z=\lbrace \delta \zeta : \zeta \in k'_{\ell'}+\Gamma^z \rbrace,
 \end{equation*}
 with   $T^{k'}_\delta$ denote the same previous translation.
\subsection{Microscopic Bidomain Model}
A vast literature exists on the bidomain modeling of the heart, we refer to \cite{char}, \cite{colli02,colli05,colli12}, \cite{henri} for more details. \\
\paragraph{\textbf{Basic equations.}} The basic equations modeling the electrical activity of the heart can be obtained as follows. First, we know that the structure of the cardiac tissue can be viewed as composed by two volumes: the intracellular space  $\Omega_i$ (inside the cells) and the extracellular space $\Omega_e$ (outside) separated by the active membrane $\Gamma^{y}$.
  
  Thus, the membrane $\Gamma^{y}$ is pierced by proteins whose role is to ensure ionic transport between the two media (intracellular and extracellular) through this membrane. So, this transport creates an electric current.\\ So by using Ohm's law, the intracellular and extracellular electrical potentials $u_{j}: \Omega_{j,T} \mapsto \R$ are related to the current volume densities $J_{j}: \Omega_{j,T} \mapsto \R^d$ for $j=i,e$ :
\begin{equation*}
J_{j}=\mathrm{M}_{j}\nabla u_{j}, \ \text{in} \ \Omega_{j,T}:=(0,T)\times\Omega_{j},
\end{equation*}
with $\mathrm{M}_{j}$ represent the corresponding conductivities of the tissue (which are assumed to be isotropic at the microscale) and are given in mS/cm$^{2}$.\\
In addition, the \textit{transmembrane} potential $v$ is known as the  potential at the membrane $\Gamma^{y}$   which is defined as follows:
\begin{equation*}
v=(u_{i}-u_{e})_{\vert_{\Gamma^{y}}} : (0,T)\times \Gamma^{y} \mapsto \R.
\end{equation*}
 
  Moreover, we assume the intracellular and extracellular spaces are source-free and thus the intracellular and extracellular potentials $u_{i}$ and $u_{e}$ are solutions to the elliptic equations:
\begin{equation}
\begin{aligned}
&-\text{div}J_{j}=0, \ \text{in} \ \Omega_{j,T}.
\end{aligned}
\label{pb}
\end{equation}
 
 According to the  current conservation law, the surface current density $\I_{m}$ is now introduced:
\begin{equation}
\I_m=-J_{i}\cdot n_i=J_{e}\cdot n_e, \ \text{on} \ \Gamma^{y}_{T}:=(0,T)\times \Gamma^{y},
\label{cond_bord} 
\end{equation}
with $n_i$ denotes the unit exterior normal to the boundary $\Gamma^{y}$ from intracellular to extracellular space and $n_e=-n_i$.

 The membrane has both a capacitive property schematized by a capacitor and a resistive property schematized by a resistor. On the one hand, the capacitive property depends on the formation of the membrane which can be represented by a capacitor of capacitance  $C_m$ (the capacity per unit area of the membrane is given in $\mu$F$/$cm$^2$). We recall that the quantity of the charge of a capacitor is $q=C_m v.$ Then, the capacitive current $\I_c$ is the amount of charge that flows per unit of time:
\begin{equation*}
\I_{c}=\pt_t q=C_m\pt_t v.
\end{equation*}
On the other hand, the resistive property depends on the ionic transport between the intracellular and extracellular media. Then, the resistive current $\I_r$ is defined by the ionic current $\I_{ion}$ measured from the intracellular  to the extracellular medium which depends on the transmembrane potential $v$ and the gating variable $w : \Gamma^{y} \mapsto \R$. 
 Since the electric current can be blocked by the membrane or can be pass through the membrane with ionic current $\I_r-\I_{app}.$ So, the charge conservation states that the total transmembrane current $\I_m$ (see \cite{colli12}) is given as follows:
\begin{equation*}
\I_m=\I_c+\I_r-\I_{app} \text{ on } \Gamma^{y}_{T},
\end{equation*} 
where $\I_{app}$ is the applied current per unit area of the membrane surface (given in $\mu$A/cm$^{2}$). 
Consequently, the transmembrane potential $v$ satisfies the following dynamic condition on $\Gamma^{y}$ involving the gating variable $w$:
\begin{equation}
\begin{aligned}
&\I_m=  C_m\pt_t v+\I_{ion}(v,w)-\I_{app}  &\text{ on } \Gamma^{y}_{T},
\\ &\pt_t w-H(v,w)=0 &\text{ on } \Gamma^{y}_{T}.
\label{cond_dyn}
\end{aligned}
\end{equation}
Herein, the functions $H$ and $\I_{ion}$ correspond to the ionic model of membrane dynamics. All surface current densities $\I_m$ and $\I_{ion}$ are given in $\mu$A/cm$^{2}$. Moreover, time is given in ms and length is given in cm.

  Mitochondria are a subcompartment of the cell bound by a
double membrane. Although some mitochondria probably do look like the traditional cigar shaped structures, it is more accurate to think of them as a budding and fusing network similar to the endoplasmic reticulum. Mitochondria are intimately involved in cellular homeostasis. Among other functions they play a part in intracellular signalling and apoptosis, intermediary metabolism, and in the metabolism of amino acids, lipids, cholesterol, steroids, and nucleotides. Perhaps most importantly, mitochondria have a fundamental role in cellular energy metabolism. This includes fatty acid oxidation, the urea cycle, and the final common pathway for ATP production-the respiratory chain (see \cite{chinnery} for more details). For this, we assume that the mitochondria are electrically insulated from the remainder of the intracellular space. Thus, we suppose that the no-flux boundary condition at the interface $\Gamma^{z}$ is given by:
\begin{equation}
\mathrm{M}_{i}\nabla u_{i}\cdot n_{z}= 0 \quad \text{ on } \Gamma^{z}_{T}:=(0,T)\times \Gamma^z,
\label{cond_dyn_z}
\end{equation}
with $n_{z}$ denotes the unit exterior normal to the boundary $\Gamma^{z}$. 


 \paragraph{\textbf{Non-dimensional analysis.}} Cardiac tissues have a number of important inhomogeneities, particularly those related to inter-cellular communications. The dimensionless analysis done correctly makes the problem simpler and clearer. In the literature, few works in that direction have been carried out, although we can cite \cite{colli12,henri,rioux} for the nondimensionalization procedure of the ionic current and \cite{richardson,whiteley} for the non-dimensional analysis in the context of bidomain equations. So, this analysis follows three steps. 
 
 First, we can define the dimensionless scale parameter:
 \begin{equation*}
 \e:=\sqrt{\dfrac{\ell^{mes}}{R_m \lambda}},
 \end{equation*}
where $R_m$ denotes the surface specific resistivity of the membrane $\Gamma^{y}$ and $\lambda:= \lambda_i+\lambda_e,$ with $\lambda_j$ represents the average eigenvalues of the corresponding conductivity  $\mathrm{M}_j$ for $j=i,e,$ over the cells' arrangement.  Now, we perform the following scaling of the space and time variables :
$$\widehat{x}=\frac{x}{L}, \quad \widehat{t}=\frac{t}{\tau}, $$
with the macroscopic units of length $L=\ell^{mes}/\e=\ell^{mic}/\delta$ and the time constant $\tau$ associated with charging the membrane by the transmembrane current is given by:
$$\tau=R_mC_m. $$
We take $\widehat{x}$ to be the variable at the macroscale (slow variable),
$$y:=\dfrac{\widehat{x}}{\e} \text{ and } z:=\dfrac{\widehat{x}}{\e\delta}$$
to be respectively the mesoscopic and microscopic space variable (fast variables) in the corresponding unit cell.
 
 Secondly, we scale all electrical potentials $u_j, \ v,$  currents and the gating variable $w$ :
 \begin{equation*}
 v=\Delta v \widehat{v}, \quad u_{j}=\Delta v\widehat{u}_{j} \text{ and }  w=\Delta w \widehat{w}_{e}
\end{equation*} 
where $\Delta v$ and $\Delta w$ are convenient units to measure electric potentials and gating variable, respectively, for $j=i,e.$ By the chain rule, we obtain:
\begin{equation*}
\dfrac{L C_m}{\tau}\pt_{\widehat{t}}\widehat{v}+\dfrac{L}{\Delta v}\left(\I_{ion}-\I_{app} \right)=-\mathrm{M}_{i}\nabla_{\widehat{x}} \widehat{u}_{i}\cdot n_i= \mathrm{M}_{e}\nabla_{\widehat{x}} \widehat{u}_{e}\cdot n_e.
\end{equation*}
Recalling that $\tau=R_m C_m$ and normalizing the conductivities $\mathrm{M}_j$ for $j=i,e$ using 
\begin{equation*}
\widehat{\mathrm{M}}_{j}=\dfrac{1}{\lambda}\mathrm{M}_{j},
\end{equation*}  
we get 
\begin{equation*}
\dfrac{L}{R_m \lambda}\pt_{\widehat{t}}\widehat{v}+\dfrac{L}{\Delta v \lambda}\left(\I_{ion}-\I_{app} \right)=-\widehat{\mathrm{M}}_{i}\nabla_{\widehat{x}} \widehat{u}_{i}\cdot n_i= \widehat{\mathrm{M}}_{e}\nabla_{\widehat{x}} \widehat{u}_{e}\cdot n_e.
\end{equation*}
Regarding the ionic functions $\I_{ion},$ $H,$ and the applied current  $\I_{app},$ we nondimensionalize them by using the following scales
\begin{equation*}
\widehat{\I}_{ion}(\widehat{v},\widehat{w})= \dfrac{R_m}{\Delta v}\widehat{\I}_{ion}(\widehat{v},\widehat{w}), \ \widehat{\I}_{app}=\dfrac{R_m}{\Delta v} \I_{app} \text{ and } \widehat{H}(\widehat{v},\widehat{w})=\dfrac{\tau}{\Delta w} H(v,w).
\end{equation*}
Consequently, we have 
\begin{equation*}
\dfrac{L}{R_m \lambda}\left( \pt_{\widehat{t}}\widehat{v}+\widehat{\I}_{ion}(\widehat{v},\widehat{w})-\widehat{\I}_{app}\right) =-\widehat{\mathrm{M}}_{i}\nabla_{\widehat{x}} \widehat{u}_{i}\cdot n_i= \widehat{\mathrm{M}}_{e}\nabla_{\widehat{x}} \widehat{u}_{e}\cdot n_e.
\end{equation*}
\begin{rem}
Recalling that the dimensionless parameter $\e,$ given by $\e:=\sqrt{\dfrac{\ell^{mes}}{R_m \lambda}},$ is the ratio between the mesoscopic cell length $\ell^{mes}$ and the macroscopic length $L$, i.e. $\e=\ell^{mes}/L$  and solving for $\e,$ we obtain
\begin{equation*}
\e=\dfrac{L}{R_m \lambda}.
\end{equation*}
\end{rem}
 
 Finally, we can convert the above microscopic bidomain system \eqref{pb}-\eqref{cond_dyn_z} to the following non-dimensional form: 
 \begin{subequations}
\begin{align}
-\nabla_{\widehat{x}}\cdot\left( \widehat{\mathrm{M}}_{i}^{\e,\delta}\nabla_{\widehat{x}} \widehat{u}_{i}^{\e,\delta}\right)  &=0 &\text{ in } \Omega_{i,T}^{\e,\delta}:=(0,T)\times\Omega_{i}^{\e,\delta}, 
\label{bid_intra}
\\ -\nabla_{\widehat{x}}\cdot\left( \widehat{\mathrm{M}}_{e}^{\e}\nabla_{\widehat{x}} \widehat{u}_{e}^{\e} \right)  &=0 &\text{ in } \Omega_{e,T}^{\e}:=(0,T)\times\Omega_{e}^\e,
\label{bid_extra}
\\ \e\left( \pt_{\widehat{t}} \widehat{v}_\e+\widehat{\I}_{ion}(\widehat{v}_\e,\widehat{w}_\e)-\widehat{\I}_{app,\e}\right) &=\widehat{\I}_m &\ \text{on} \ \Gamma_{\e,T}:=(0,T)\times\Gamma_{\e},
\label{bid_onface}
\\-\widehat{\mathrm{M}}_{i}^{\e,\delta}\nabla_{\widehat{x}} \widehat{u}_{i}^{\e,\delta} \cdot n_i=\widehat{\mathrm{M}}_{e}^{\e}\nabla_{\widehat{x}} \widehat{u}_{e}^{\e} \cdot n_e & =\widehat{\I}_m &\ \text{on} \ \Gamma_{\e,T},
\label{bid_mesocont}
\\ \pt_{\widehat{t}} \widehat{w}_\e-\widehat{H}(\widehat{v}_\e,\widehat{w}_\e) &=0 & \text{ on }  \Gamma_{\e,T},
\label{bid_dyn}
\\ \widehat{\mathrm{M}}_{i}^{\e,\delta}\nabla_{\widehat{x}} \widehat{u}_{i}^{\e,\delta} \cdot n_{z} & =0 &\ \text{on} \ \Gamma_{\delta,T},
\label{bid_microcont}
\end{align}
\label{pbscale}
\end{subequations}
with each equation correspond to the following sense:
\eqref{bid_intra} Intra quasi-stationary conduction, \eqref{bid_extra} Extra quasi-stationary conduction, \eqref{bid_onface} Reaction onface condition, \eqref{bid_mesocont} Meso-continuity equation, \eqref{bid_dyn} Dynamic coupling, \eqref{bid_microcont} Micro-boundary condition.

For convenience, the superscript $ \  \widehat{\cdot} \ $ of the dimensionless variables is omitted. 
Note that the bidomain equations are invariant with respect to the scaling parameters $\e$ and $\delta$. Then, we define the rescaled electrical potential as follows:
$$u_{i}^{\e,\delta}(t,x):=u_{i}\left( t, x,\frac{x}{\e},\frac{x}{\e\delta}\right), \ \ u_{e}^{\e}(t,x):= u_{e}\left( t, x,\frac{x}{\e}\right). $$
Analogously, we obtain the rescaled transmembrane potential  $v_\e=(u_{i}^{\e,\delta}-u_{e}^{\e}){\vert_{\Gamma_{\e,T}}}$ and gating variable $w_\e.$ Thus, we define also the following rescaled conductivity matrices: 
\begin{equation}
\mathrm{M}_{i}^{\e,\delta}(x):=\mathrm{M}_{i}\left( x,\frac{x}{\e},\frac{x}{\e\delta}\right) \ \text{ and } \ \mathrm{M}_{e}^{\e}(x):= \mathrm{M}_{e}\left( x,\frac{x}{\e}\right)
\label{M_ie}
\end{equation}
  satisfying the elliptic and periodicity conditions \eqref{cond_Me}-\eqref{cond_Mi}. 
  
    Finally, The ionic current $\I_{ion}(v,w)$ can be decomposed into $\mathrm{I}_{1,ion}(v): \R \rightarrow \R$ and $\mathrm{I}_{2,ion}(w); \R \rightarrow \R$, where $\I_{ion}(v,w)=\mathrm{I}_{1,ion}(v)+\mathrm{I}_{2,ion}(w)$. Furthermore, $\mathrm{I}_{1,ion}$ is considered as a $C^1$ function, $\mathrm{I}_{2,ion}$ and $H : \R^2 \rightarrow \R$ are linear functions. Also, we assume that there exists $r\in (2,+\infty)$ and constants $\alpha_1,\alpha_2,\alpha_3, \alpha_4,\alpha_5, C>0$ and $\beta_1, \beta_2>0$ such that:
\begin{subequations}
\begin{align}
&\dfrac{1}{\alpha_1} \abs{v}^{r-1}\leq \abs{\mathrm{I}_{1,ion}(v)}\leq \alpha_1\left(  \abs{v}^{r-1}+1\right), \,\abs{\mathrm{I}_{2,ion}(w)}\leq \alpha_2(\abs{w}+1), 
\\& \abs{H(v,w)}\leq \alpha_3(\abs{v}+\abs{w}+1),\text{ and }
    \mathrm{I}_{2,ion}(w)v-\alpha_4H(v,w)w\geq \alpha_5 \abs{w}^2,
\\& \tilde{\mathrm{I}}_{1,ion} : z\mapsto \mathrm{I}_{1,ion}(z)+\beta_1 z+\beta_2 \text{ is strictly increasing with } \lim \limits_{z\rightarrow 0} \tilde{\mathrm{I}}_{1,ion}(z)/z=0,
\\& \forall z_1,z_2 \in \R,\,\,\left(\tilde{\mathrm{I}}_{1,ion}(z_1)-\tilde{\mathrm{I}}_{1,ion}(z_2) \right)(z_1-z_2)\geq \dfrac{1}{C} \left(1+\abs{z_1}+\abs{z_2} \right)^{r-2} \abs{z_1-z_2}^{2}. 
\end{align}
\label{A_H_I_asym}
\end{subequations}
\begin{rem} 
In the mathematical analysis of bidomain equations, several paths have been followed in the literature according to the definition of the ionic currents. We summarize below the encountered various cases:
\begin{enumerate}\item[$1.$] \textbf{Physiological models}\\
These types of models attempt to describe specific actions within the cell membrane. Such exact models are derived either by fitting the parameters of an  equation to match experimental data or by defining equations that were confirmed by later experiments. Moreover, they are based on the cell membrane formulation developed by Hodgkin and Huxley for nerve fibers \cite{hodgkin} (see \cite{noble12} for more details). To go further in the physiological description, some models consider the concentrations as variables of the system, see for example the Beeler-Reuter model \cite{beeler} and the Luo-Rudy model \cite{luo1,luo2,luo3}. In \cite{marco06,marco}, such models are considered.
\item[]
\item[$2.$] \textbf{Phenomenological models}\\
Other non-physiological models have been introduced as approximations of ion current models. They can be used in large problems because they are typically small and fast to solve, although they are less flexible in their response to variations in cellular properties such as concentrations or cell size. We take in this paper the FitzHugh-Nagumo model \cite{fitz,nagumo} that satisfies assumptions \eqref{A_H_I_asym} which reads as
  \begin{subequations}
  \begin{align}
  H(v,w)&= av-bw, \\  \I_{ion}(v,w)&=\left( \lambda v (1-v)(v-\theta) \right) + (-\lambda w):= \mathrm{I}_{1,ion}(v)+\mathrm{I}_{2,ion}(w)
  \end{align}
  \label{ionic_model_asym}
  \end{subequations}
  where $a, b, \lambda, \theta$ are given parameters with $a,b\geq 0, \ \lambda<0$ and $0<\theta<1.$ According to this model, the functions $\I_{ion}$ and $H$ are continuous and the non-linearity $\mathrm{I}_{a,ion}$ is of cubic growth at infinity then the most appropriate value is $r=4.$ We end this remark by mentioning other reduced ionic models: the Roger-McCulloch model \cite{roger} and the Aliev-Panfilov model \cite{aliev}, may consider more general that the previous model but still rise some mathematical difficulties.
\end{enumerate}
\end{rem}
We complete system \eqref{pbscale} with no-flux boundary conditions: 
\begin{equation*}
\left( \mathrm{M}_{i}^{\e,\delta}\nabla u^{\e,\delta}_{i}\right) \cdot \mathbf{n}=\left( \mathrm{M}_{e}^{\e}\nabla u_{e}^{\e}\right) \cdot\mathbf{n}=0 \ \text{ on } \  (0,T)\times \pt_{\text{ext}} \Omega,
\end{equation*}
and appropriate the initial Cauchy conditions for transmembrane potential $v$ and gating variable $w$. Herein, $\mathbf{n}$ is the outward unit normal to the exterior boundary of $\Omega.$
 
  Clearly, the equations in \eqref{pbscale} are invariant under the simultaneous change of $u_i^{\e,\delta}$ and $u_{e}^{\e}$ into $u_i^{\e,\delta}+k;$ $u_{e}^{\e}+k,$ for any $k\in \R.$ Hence, we may impose the following normalization condition:
\begin{equation}
\int_{\Omega_{e}^{\e}}u_{e}^{\e}(t,x)dx=0 \text{ for a.e. } t\in(0,T).
\label{normalization_cond} 
\end{equation}
\section{Asymptotic Expansion Homogenization}\label{asym}
In this section, we will introduce a homogenization method based on asymptotic expansion using multi-scale variables (i.e. slow and fast variables). The aim is to show how to obtain a mathematical writing of the macroscopic model from the microscopic model. This method, among others, is a formal and intuitive method for predicting the mathematical writing of a homogenized solution that can eventually approach the solution of the initial problem \eqref{pbscale}. 
 
 For that, we start to treat the problem in the extracellular medium then we will solve the other one in the intracellular medium using this method.
 
 \subsection{Extracellular problem}\label{asymextra}In the literature, Cioranescu and Donato \cite{doina} are applied and developed the two-scale asymptotic expansion method established by Benssousan and Papanicolaou \cite{ben} on a problem defined at two scales  to obtain the homogenized model (see also \cite{sanchez,olei,bakhvalov}). Further, the authors in \cite{burridge} have been used this method to derive the macroscopic linear behavior of a porous elastic solid saturated with a compressible viscous fluid. Its derivation is based on the linear elasticity equations in the solid, the linearized Navier-Stokes equations in the fluid, and the appropriate conditions at the solid-fluid boundary.
 
  In our approach, we investigate the same two-scale technique for the extracellular problem. Whereas for the intracellular domain, we develop a three-scale approach applied to the intracellular problem to handle with the two structural levels of this domain (see Figure \ref{three_scale}).  We recall the following initial extracellular problem:
\begin{equation}
 \begin{aligned}
\A_\e u_{e}^\e &=0 &\ \text{ in } \Omega_{e,T}^{\e}, \\ \mathrm{M}_{e}^{\e}\nabla u_{e}^{\e} \cdot n_e=\e\left( \pt_{t} v_\e+\I_{ion}(v_\e,{w}_\e)-\I_{app,\e}\right) &=\I_m &\ \text{on} \ \Gamma_{\e,T}, 
 \end{aligned}
 \label{pbiniextra}
 \end{equation}
 with $\A_\e=-\nabla \cdot\left( \mathrm{M}_{e}^{\e}\nabla\right),$ where the extracellular conductivity matrices $\mathrm{M}^{\e,\delta}_{e}$ are defined by: $$\mathrm{M}^{\e}_{e}(x)=\mathrm{M}_{e}\left(\dfrac{x}{\e}\right), \
  a.e. \ \text{on} \ \R^d,$$

 satisfying the following elliptic and periodic conditions: 
 \begin{equation}
 \begin{cases}
\mathrm{M}_e(y)\in M(\alpha,\beta,Y),\\ \mathrm{M}_e=(\mathrm{m}^{pq}_{e})_{1 \leq p,q \leq d} \text{ with }  \mathrm{m}^{pq}_{e} \ y\text{-periodic}, \ \forall p,q=1,\dots,d,  
 \end{cases}
 \label{cond_Me}
\end{equation}
with $\alpha, \beta \in \R,$ such that $0<\alpha<\beta$ and $M(\alpha,\beta,Y)$ given by  Definition \ref{M}.
 
The two-scale asymptotic expansion is assumed for the electrical potential $u^\e_e$ as follows:
\begin{equation}
 u^\e_e (t,x):=u_e\left( t, x,\dfrac{x}{\e}\right) =u_{e,0}\left(t,x,\dfrac{x}{\e}\right)+\e u_{e,1} \left(t,x,\dfrac{x}{\e}\right)+\e^{2} u_{e,2} \left(t,x,\dfrac{x}{\e}\right)+ \cdots
 \label{ue_extra}
 \end{equation}
 where each $u_{j}(\cdot,y)$  is $y$-periodic function dependent on time $t\in(0,T),$ slow (macroscopic) variable  $x$ and the fast (mesoscopic) variable $y$. The slow and fast variables correspond respectively to the global and local structure of the field. Similarly, the applied current $\I_{app,\e}$ has the same two-scale asymptotic expansion.\\
Consequently, the full operator $\A_\e$ in the initial problem  \eqref{pbiniextra} is represented as:
\begin{equation}
\A_\e u^\e_e(t,x)=[(\e^{-2} \A_{yy}+\e^{-1} \A_{xy}+\e^{0} \A_{xx})u_e]\left( t,x,\dfrac{x}{\e} \right)
\label{Aeue_extra}
\end{equation}
with each operator defined by: 
\begin{equation*}
 \begin{cases}\A_{yy}=-\overset{d}{\underset{p,q=1}{\sum}}\dfrac{\pt}{\pt y_p}\left(\mathrm{m}^{pq}_{e}(y)\dfrac{\pt}{\pt y_q}\right), \\ \A_{xy}=-\overset{d}{\underset{p,q=1}{\sum}}\dfrac{\pt}{\pt y_p}\left(\mathrm{m}^{pq}_{e}(y)\dfrac{\pt}{\pt x_q}\right)-\overset{d}{\underset{p,q=1}{\sum}}\dfrac{\pt}{\pt x_p}\left(\mathrm{m}^{pq}_{e}(y)\dfrac{\pt}{\pt y_q}\right),\\\A_{xx}=-\overset{d}{\underset{p,q=1}{\sum}}\dfrac{\pt}{\pt x_p}\left(\mathrm{m}^{pq}_{e}(y)\dfrac{\pt}{\pt x_q}\right).
 \end{cases}
\end{equation*}

Now, we substitute the asymptotic expansion \eqref{ue_extra} of $u^\e_e$ in the developed operator \eqref{Aeue_extra} to obtain
\begin{equation*}
\begin{aligned}
\A_\e u^\e_e(x) &=[\e^{-2} \A_{yy} u_{e,0}+\e^{-1} \A_{yy} u_{e,1}+\e^{0} \A_{yy} u_{e,2}+ \cdots]\left( t, x,\dfrac{x}{\e}\right) \\ & \ \ +[\e^{-1} \A_{xy} u_{e,0}+\e^{0} \A_{xy} u_{e,1}+ \cdots]\left( t, x,\dfrac{x}{\e}\right)\\ & \ \ +[\e^{0} \A_{xx} u_{e,0}+ \cdots]\left( t, x,\dfrac{x}{\e}\right)
\\&=[\e^{-2}\A_{yy} u_{e,0}+\e^{-1}\left(\A_{yy} u_{e,1}+\A_{xy} u_{e,0} \right)
\\ & \ \ +\e^{0}\left(\A_{yy} u_{e,2}+\A_{xy} u_{e,1}+\A_{xx} u_{e,0} \right)]\left( t, x,\dfrac{x}{\e} \right)+\cdots.
\end{aligned}
\end{equation*}

Similarly,  we substitute the asymptotic expansion \eqref{ue_extra} of $u^\e_e$ into the boundary condition equation \eqref{pbiniextra} on $\Gamma^{y}$. Consequently, by equating the powers-like terms of $\e^\ell$ to zero ($\ell=-2,-1,0$), we have to solve the following system of equations for the functions $u_{e,k}(t,x,y), \ k=0,1,2:$

 \begin{equation}
 \begin{cases}
\A_{yy} u_{e,0}=0 \ \text{in} \ Y_{e},\\ u_{e,0} \ y\text{-periodic}, \\  \mathrm{M}_{e}\nabla_y u_{e,0}  \cdot n_e= 0 \ \text{on} \ \Gamma^{y},  
 \end{cases}
 \label{Ayyue0}
 \end{equation}

 \begin{equation}
 \begin{cases}
\A_{yy} u_{e,1}=-\A_{xy} u_{e,0} \ \text{in} \ Y_{e},\\ u_{e,1} \ y\text{-periodic}, \\ \left( \mathrm{M}_{e}\nabla_y u_{e,1}+ \mathrm{M}_{e}\nabla_x u_{e,0}\right)  \cdot n_e= 0 \ \text{on} \ \Gamma^{y}, 
 \end{cases}
 \label{Ayyue1}
 \end{equation}
 
 \begin{equation}
 \begin{cases}
\A_{yy} u_{e,2}=-\A_{xy} u_{e,1}-\A_{xx} u_{e,0} \ \text{in} \ Y_{e},\\ u_{e,2} \ y\text{-periodic}, \\ \left( \mathrm{M}_{e}\nabla_y u_{e,2}+ \mathrm{M}_{e}\nabla_x u_{e,1}\right)  \cdot n_e= \pt_{t} v_0+\I_{ion}(v_0,{w}_0)-\I_{app,0} \ \text{on} \ \Gamma^{y},
 \end{cases}
 \label{Ayyue2}
 \end{equation}
 
The authors in \cite{ben}-\cite{doina} have successively solved the three systems into Dirichlet boundary conditions \eqref{Ayyue0}-\eqref{Ayyue2}. Herein, the functions $u_{e,0}, u_{e,1}$ and $u_{e,2}$ in the asymptotic expansion \eqref{ue_extra} for the extracellular potential $u^\e_{e}$ satisfy the Neumann boundary value problems \eqref{Ayyue0}-\eqref{Ayyue2} in the local portion $Y_e$ of a unit cell $Y$ (see \cite{colli02,henri} for the case of Laplace equations).\\ The resolution is described as follows:  
\begin{enumerate}
\item[$\bullet$ \textbf{First step} ]We begin with the first boundary value problem \eqref{Ayyue0} whose variational formulation: 
\begin{equation}
 \begin{cases}
\text{Find} \ \dot{u}_{e,0} \in \W_{per}(Y_{e}) \ \text{such that} \\ \dot{a}_{Y_{e}}(\dot{u}_{e,0},\dot{v})= \displaystyle \int_{\pt Y_{e}} (\mathrm{M}_e \nabla_y u_{e,0} \cdot n_e) v \ d\sigma_y,\ \forall \dot{v}\in \W_{per}(Y_{e}),
 \end{cases}
 \label{FvAyyue0}
\end{equation} 
with $\dot{a}_ {Y_{e}}(\dot{u},\dot{v})$ given by: 
\begin{equation}
\dot{a}_ {Y_{e}}(\dot{u},\dot{v})=\int_{Y_{e}} \mathrm{M}_e \nabla_y u \nabla_y v dy, \ \forall u \in\dot{u}, \ \forall v \in\dot{v},\ \forall \dot{u},\ \forall \dot{v} \in \W_{per}(Y_{e})
\label{dotaye}
\end{equation}
 and $\W_{per}(Y_{e})$ is given by Definition \ref{w_per}.\\ We want to clarify the right hand side of the variational formulation \eqref{FvAyyue0}. By the definition of $\pt Y_{e}:=(\pt_{\text{ext}} Y  \cap \, \pt Y_{e})\, \cup \, \Gamma^{y},$ we use Proposition \ref{cond_Y_per} and the $y$-periodicity of $\mathrm{M}_i$ by taking account the boundary condition on $\Gamma^{y}$ to say that :
   \begin{align*}
\int_{\pt Y_{e}} &(\mathrm{M}_{e} \nabla_y u_{e,0} \cdot n_e) v \ d\sigma_y \\&=\int_{\pt_{\text{ext}} Y  \cap \pt Y_{e}} (\mathrm{M}_{e} \nabla_y u_{e,0} \cdot n_e) v \ d\sigma_y+\int_{\Gamma^{y}} (\mathrm{M}_{e} \nabla_y u_{e,0} \cdot n_e) v \ d\sigma_y =0.
 \end{align*}
 Using Theorem \ref{cond_per_neu}, we can prove the existence and uniqueness of the solution $\dot{u}_{e,0}$. Then, the problem \eqref{Ayyue0} has a unique solution $u_{e,0}$ independent of $y,$ so  we deduce that: $$u_{e,0}(t,x,y)=u_{e,0}(t,x).$$
In the next section, we show that $u_{i,0}$ does not depend on $y$ and $z$ (by the same strategy). Since $v_0=(u_{i,0}-u_{e,0})\vert_{\Gamma^{y}}$ then we also deduce that $v_0$ and $w_0$ not depend on the mesoscopic variable $y.$ 

\item[$\bullet$ \textbf{Second step}] We now turn to the second boundary value problem \eqref{Ayyue1}. 
Using Theorem \ref{cond_per_neu}, we obtain that the second system \eqref{Ayyue1} 
has a unique weak solution $\dot{u}_{e,1} \in \W_{per}(Y_{e})$ (defined by \cite{ben} and \cite{olei}). 

Thus, the linearity of terms in the right hand side of equation \eqref{Ayyue1} 
 suggests to look for $\dot{u}_{e,1}$ under the following form:
\begin{equation}
\dot{u}_{e,1}(t,x,y)=\overset{d}{\underset{q=1}{\sum}}\dot{\chi}_e^q(y)\dfrac{\pt \dot{u}_{e,0}}{\pt x_q}(t,x) \ \text{in} \ \W_{per}(Y_{e}),
\label{dotue1}
\end{equation}
with the corrector function $\dot{\chi}_e^q$ satisfies the following $\e$-cell problem:
\begin{equation}
\begin{cases}
\A_{yy} \dot{\chi}_e^q=\overset{d}{\underset{p=1}{\sum}}\dfrac{\pt \mathrm{m}^{pq}_{e}}{\pt y_p} \ \text{in} \ Y_{e},\\ \dot{\chi}_e^q \ y\text{-periodic}, \\ \mathrm{M}_{e} \nabla_y \dot{\chi}_e^q \cdot n_e= - (\mathrm{M}_{e}e_q )\cdot n_e \text{ on } \Gamma^{y}, 
 \end{cases}
 \label{Ayychie}
 \end{equation}
 for $e_q$, $q=1,\dots,d,$ the standard canonical basis in $\R^d.$ Moreover, we can choose a representative element $\chi_e^q$ of the class $\dot{\chi}_e^q$  satisfying the following variational formulation:
 \begin{equation}
 \begin{cases}
\text{Find} \ \chi_e^q \in W_{\#}(Y_{e}) \ \text{such that} \\ \\ a_{Y_{e}}(\chi_e^q,v)=(F,v)_{(W_{\#}(Y_{e}))',W_{\#}(Y_{e})}, \ \ \forall v \in W_{\#}(Y_{e}),
 \end{cases}
 \label{FvAyychie}
\end{equation}
with $a_{Y_{e}}$ given by \eqref{dotaye} and $F$ defined by: 
$$(F,v)_{(W_{\#}(Y_{e}))',W_{\#}(Y_{e})}=\overset{d}{\underset{p=1}{\sum}}\displaystyle\int_{Y_{e}} \mathrm{m}^{pq}_{e}(y) \dfrac{\pt v}{\pt y_p}dy,$$
where the space $W_{\#}(Y_e)$ is given by the expression \eqref{W}. Since $F$ belongs to $(W_{\#}(Y_{e}))'$ then the condition of Theorem \ref{cond_per_neu} is imposed in order to guarantee existence and uniqueness of the solution. \\Thus, by the form of $\dot{u}_{e,1}$  given by \eqref{dotue1}, the solution $u_{e,1}$ of the second system \eqref{Ayyue1} can be represented by the following ansatz:
\begin{equation} 
u_{e,1}(t,x,y)=\chi_e(y)\cdot \nabla_x u_0(t,x)+\tilde{u}_{e,1}(t,x) \ \text{with} \ u_{e,1} \in \dot{u}_{e,1},
\label{ue1}
\end{equation}
where $\tilde{u}_{e,1}$ is a constant with respect to $y$  (i.e $ \tilde{u}_{e,1} \in \dot{0}  $ in $\W_{per}(Y)).$
\item[$\bullet$ \textbf{Last step} ] We now pass to the last boundary value problem \eqref{Ayyue2}. Taking into account the form of $u_{e,0}$ and $u_{e,1}$,  we obtain
\begin{equation*}
\begin{aligned}
&-\A_{xy} u_{e,1}-\A_{xx} u_{e,0} 
\\& \qquad=\overset{d}{\underset{p,q=1}{\sum}}\dfrac{\pt}{\pt y_p}\left(\mathrm{m}^{pq}_{e}(y)\dfrac{\pt u_{e,1}}{\pt x_q}\right)+\overset{d}{\underset{p,q=1}{\sum}}\dfrac{\pt}{\pt x_p}\left(\mathrm{m}^{pq}_{e}(y)\left(\dfrac{\pt u_{e,1}}{\pt y_q}+\dfrac{\pt u_{e,0}}{\pt x_q}\right)\right). 
\end{aligned}
\end{equation*}
Consequently, this system \eqref{Ayyue2} have the following variational formulation:
\begin{equation}
 \begin{cases}
\text{Find} \ \dot{u}_{e,2} \in \W_{per}(Y_{e}) \ \text{such that} \\ \dot{a}_{Y_{e}}(\dot{u}_{e,2},\dot{v})=(F_2,\dot{v})_{(\W_{per}(Y_{e}))',\W_{per}(Y_{e})} \ \forall \dot{v}\in \W_{per}(Y_{e}),
 \end{cases}
 \label{FvAyyue2}
\end{equation}
with $\dot{a}_{Y_{e}}$ given by \eqref{dotaye} and $F_2$ defined by:
\begin{equation}
\begin{aligned}
&(F_2,\dot{v})_{(\W_{per}(Y_e))',\W_{per}(Y_e)}
\\& \ \ \ =\int_{\Gamma^{y}} \left( \mathrm{M}_{e}\nabla_y u_{e,2}+ \mathrm{M}_{e}\nabla_x u_{e,1}\right)  \cdot n_e \ v \ d\sigma_y -\overset{d}{\underset{p,q=1}{\sum}}\int_{Y_e} \mathrm{m}^{pq}_{e}(y)\dfrac{\pt u_{e,1}}{\pt x_q}\dfrac{\pt v}{\pt y_p} dy
\\ & \quad \ \ +\overset{d}{\underset{p,q=1}{\sum}} \int_{Y_e} \dfrac{\pt}{\pt x_p}\left(\mathrm{m}^{pq}_{e}(y)\left(\dfrac{\pt u_{e,1}}{\pt y_q}+\dfrac{\pt u_{e,0}}{\pt x_q}\right) \right) v dy,  \ \ \forall v \in\dot{v}, \ \forall \dot{v} \in \W_{per}(y).
\end{aligned}
\label{F2e}
\end{equation}
The problem \eqref{FvAyyue2}-\eqref{F2e} is well-posed according to Theorem \ref{cond_per_neu} under the compatibility condition: 
\begin{equation*}
(F_2,1)_{(\W_{per}(Y_e))',\W_{per}(Y_e)}=0.
\end{equation*}
which equivalent to:
\begin{equation*}
-\overset{d}{\underset{p,q=1}{\sum}} \int_{Y_e} \dfrac{\pt}{\pt x_p}\left(\mathrm{m}^{pq}_{e}(y)\left(\dfrac{\pt u_{e,1}}{\pt y_q}+\dfrac{\pt u_{e,0}}{\pt x_q}\right) \right)dy=\abs{\Gamma^{y}} \left( \pt_{t} v_0+\I_{ion}(v_0,{w}_0)-\I_{app}\right).
\end{equation*}
In addition, we replace $u_{e,1}$ by its form \eqref{ue1} in the above condition to obtain:
\begin{align*}
&-\overset{d}{\underset{p,q=1}{\sum}} \int_{Y_e} \dfrac{\pt}{\pt x_p}\left(\mathrm{m}^{pq}_{e}(y)\left(\overset{d}{\underset{k=1}{\sum}}\dfrac{\pt \chi_e^k}{\pt y_q}\dfrac{\pt u_{e,0}}{\pt x_k}+\dfrac{\pt u_{e,0}}{\pt x_q}\right) \right)dy
\\ & \qquad \qquad=\abs{\Gamma^{y}} \left( \pt_{t} v_0+\I_{ion}(v_0,{w}_0)-\I_{app}\right).
\end{align*} 
By expanding the sum and permuting the index, we obtain
\begin{align*}
&-\overset{d}{\underset{p,q=1}{\sum}}\overset{d}{\underset{k=1}{\sum}} \int_{Y_e} \dfrac{\pt}{\pt x_p}\left(\mathrm{m}^{pq}_{e}(y)\dfrac{\pt \chi_e^k}{\pt y_q}\dfrac{\pt u_{e,0}}{\pt x_k} \right)dy-\overset{d}{\underset{p,k=1}{\sum}} \int_{Y_e} \dfrac{\pt}{\pt x_p}\left(\mathrm{m}_e^{pk}(y)\dfrac{\pt u_{e,0}}{\pt x_k} \right)dy
\\ & \qquad\qquad=\abs{\Gamma^{y}} \left( \pt_{t} v_0+\I_{ion}(v_0,{w}_0)-\I_{app}\right)
\end{align*}
which equivalent to find $u_{e,0}$ satisfying the following problem:
\begin{align*}
 &-\overset{d}{\underset{p,k=1}{\sum}}\left[\dfrac{1}{\abs{Y}}\overset{d}{\underset{q=1}{\sum}} \int_{Y_e} \left( \mathrm{m}_e^{pk}(y)+\mathrm{m}^{pq}_{e}(y)\dfrac{\pt \chi_e^k}{\pt y_q}\right)  dy\right] \dfrac{\pt^2 u_{e,0}}{\pt x_p \pt x_k}
 \\& \qquad \qquad=\dfrac{\abs{\Gamma^{y}}}{\abs{Y}}\left(  \pt_{t} v_0+\I_{ion}(v_0,{w}_0)-\I_{app}\right),
\end{align*}
where $\I_{app}(t,x)=\dfrac{1}{\abs{\Gamma^{y}}}\displaystyle\int_{\Gamma^{y}}\I_{app,0}(\cdot,y) \ d\sigma_y.$ 
\end{enumerate}
 Consequently, we see that's exactly the \textbf{homogenized} equation satisfied by $u_{e,0}$ of the extracellular problem can be rewritten as:
\begin{equation}
\B_{xx} u_{e,0}= \mu_{m} \left( \pt_{t} v_0+\I_{ion}(v_0,{w}_0)-\I_{app}\right)  \ \text{on} \ \Omega_T,
 \label{pbhomex}
 \end{equation}
 where $\mu_{m}=\abs{\Gamma^{y}}/\abs{Y}.$
Herein, the homogenized operator $\B_{xx}$ is defined by :
 \begin{equation}
 \B_{xx}=-\nabla_x \cdot \left(  \widetilde{\mathbf{M}}_e \nabla_x \right)=-\overset{d}{\underset{p,k=1}{\sum}}\dfrac{\pt}{\pt x_p}\left(\widetilde{\mathbf{m}}^{pq}_{e}\dfrac{\pt}{\pt x_k}\right)
 \label{Bxxhome}
 \end{equation}
 with the coefficients of the homogenized conductivity matrices $\widetilde{\mathbf{M}}_e=\left( \widetilde{\mathbf{m}}^{pk}_e\right)_{1\leq p,k \leq d}$ defined by:
\begin{equation}
\widetilde{\mathbf{m}}^{pk}_e:=\dfrac{1}{\abs{Y}}\overset{d}{\underset{q=1}{\sum}}\displaystyle\int_{Y_e}\left( \mathrm{m}_e^{pk}+\mathrm{m}^{pq}_{e}\dfrac{\pt \chi_e^k}{\pt y_q}\right) \ dy.
\label{Mb_e}
\end{equation} 

\begin{rem}[Comparison with other papers]
The technique we use in the extracellular problem is closely related to that of Krassowska and Neu \cite{neukra}, with some clarifications, although the resulting model differs in important ways (described in Subection \ref{asymintra}). Keener and Panfilov \cite{keener96} consider a network of myocytes, and transform to a local curvilinear coordinate system in which one coordinate is aligned with the fiber orientation. They make a transformation to the reference frame and then obtain the bidomain model analogous to that performed by Krassowska and Neu \cite{neukra} on a regular lattice of myocytes. As such, this model provides insight into the mechanism of direct stimulation and defibrillation of cardiac tissue after injection of large currents. Further, Richardson and Chapman \cite{richardson} have been applied the two-scale asymptotic expansion to bidomain problems which have an almost periodic micro-structure not in Cartesian coordinates but in a general curvilinear coordinate system. They used this method to derive a version of the bidomain equations describing the macroscopic electrical activity of cardiac tissue. The treatment systematically took into account the non-uniform orientation of the cells in the tissue and the deformation of the tissue due to the heart beat. Recently, Whiteley \cite{whiteley} used the homogenization technique for an almost periodic micro-structure described by Richardson and Chapman \cite{richardson}, to derive the tissue level bidomain equations. They also presented some observations on the entries of the conductivity tensors, as well as some observations arising from the computation of the numerical solution of $\e$-cell problems.
\end{rem}

\subsection{Intracellular problem}\label{asymintra} Using the two-scale asymptotic expansion method, the extracellular problem is treated on two scales. Our derivation bidomain model is based on a new three-scale approach. We apply a three-scale asymptotic expansion in the intracellular problem to obtain its homogenized equation. Recall that $u^{\e,\delta}_i$ the solution of the following initial intracellular problem:  
 \begin{equation}
 \begin{aligned}
\A_{\e,\delta} u_{i}^{\e,\delta} &=0 &\ \text{ in } \Omega_{i,T}^{\e,\delta}, \\ -\mathrm{M}_{i}^{\e,\delta}\nabla u_{i}^{\e,\delta} \cdot n_i=\e\left( \pt_{t} v_\e+\I_{ion}( v_\e, {w}_\e)-\I_{app,\e}\right) &=\I_m &\ \text{on} \ \Gamma_{\e,T},
\\ -\mathrm{M}_{i}^{\e,\delta}\nabla u_{i}^{\e,\delta} \cdot n_{z} & =0 &\ \text{on} \ \Gamma_{\delta,T}, 
 \end{aligned}
 \label{pbinintra}
 \end{equation}
 with $\A_{\e,\delta}=-\nabla \cdot \left( \mathrm{M}_{i}^{\e,\delta}\nabla\right),$ where the intracellular conductivity matrices $\mathrm{M}^{\e,\delta}_{i}$ are defined by: \begin{equation*}
 \mathrm{M}^{\e,\delta}_{i}(x)=\mathrm{M}_{i}\left(\dfrac{x}{\e}, \dfrac{x}{\e\delta}\right),
\end{equation*}  
 satisfying the following elliptic and periodicity conditions: 
\begin{equation}
 \begin{cases}
\mathrm{M}_i(y,\cdot) \in M(\alpha,\beta,Y),\ \mathrm{M}_i( \cdot, z) \in M(\alpha,\beta,Z),\\ \mathrm{M}_i=\left( \mathrm{m}_{i}^{pq}\right)_{1 \leq p,q \leq d} \text{ with }  \mathrm{m}^{pq}_{i} \ y \text{- and } z  \text{-periodic}, \ \forall p,q=1,\dots,d,  
 \end{cases}
 \label{cond_Mi}
\end{equation}
with $\alpha, \beta \in \R,$ such that $0<\alpha<\beta$ and $M(\alpha,\beta,\mathcal{O})$ given by Definition \ref{M}.
 
 In the intracellular problem, we consider three different scales: the slow variable $x$ describes the macroscopic one, the fast variables $\dfrac{x}{\e}$ describes the mesoscopic one while  $\dfrac{x}{\e\delta}$ describes the microscopic one.

 To proceed with multi-scale formulation of the microscopic bidomain problem, a three-scale asymptotic expansion is assumed for the intracellular potential $u_{i}^{\e,\delta}$ as follows:
\begin{equation}
\begin{aligned}
u^{\e,\delta}_i ( t, x):= u_i \left( t, x, \dfrac{x}{\e}, \dfrac{x}{\e\delta}\right) &=u_{i,0}\left( t, x,\dfrac{x}{\e}, \dfrac{x}{\e\delta}\right)+\e u_{i,1} \left( t, x,\dfrac{x}{\e},\dfrac{x}{\e \delta}\right)+\e \delta u_{i,2} \left( t, x,\dfrac{x}{\e},\dfrac{x}{\e \delta}\right)
 \\ & \ \ \ +\e^2 u_{i,3} \left( t,x,\dfrac{x}{\e},\dfrac{x}{\e \delta}\right)+\e^2 \delta u_{i,4} \left( t, x,\dfrac{x}{\e},\dfrac{x}{\e \delta} \right) 
 \\& \ \ \ +\e^2 \delta^2 \ u_{i,5} \left( t, x,\dfrac{x}{\e},\dfrac{x}{\e \delta}\right)+ \cdots
 \end{aligned}
 \label{ue_intra}
 \end{equation}
where each $u_{i,q}(\cdot,y,z)$  is $y$- and $z$-periodic function dependent on time $t\in(0,T),$ the macroscopic variable  $x,$ the  mesoscopic variable $y,$ and the microscopic variable $z.$ \\ Next, we use the chain rule to derive with respect to $x$ 
$$\dfrac{\pt u^{\e,\delta}_i}{\pt x_q}( t, x)=\left[ \dfrac{\pt u_i}{\pt x_q} +\dfrac{1}{\e}\dfrac{\pt u_i}{\pt y_q}+\dfrac{1}{\e \delta}\dfrac{\pt u_i}{\pt z_q} \right]\left(t, x, \dfrac{x}{\e}, \dfrac{x}{\e \delta} \right).$$
\begin{rem}
The authors in \cite{ramirez} used the iterated three-scale homogenization methods to study macroscopic performance of hierarchical composites in the context of mechanics where the microscale and mesoscale are very well-separated, i.e.
$$u^{\e,\delta}(x,y,z)=u_0(x,y,z)+\overset{\infty}{\underset{k=1}{\sum}}\e^{k} u_k(x,y,z)+\overset{\infty}{\underset{k=1}{\sum}}\delta^{k} u'_k(x,y,z),$$ with $y=x/\e$ and $z=x/\delta$ $(\delta<<\e)$.  The approach proposed in the present work is exploited the effective properties of cardiac tissue with multiple small-scale configurations.  We note that our present technique recovers the classical reiterated homogenization \cite{ben} where $\delta=\e$.
\end{rem}
Consequently, we can write the full operator $\A_{\e,\delta}$ in the initial problem \eqref{pbinintra} as follows:
\begin{equation}
\begin{aligned}
\A_{\e,\delta} u^{\e,\delta}_i( t, x)&=-\left[ \nabla \cdot  \left( \mathrm{M}^{\e,\delta}_{i} \nabla u^{\e,\delta}_{i}\right) \right]( t, x)  
\\&=-\left[\overset{d}{\underset{p,q=1}{\sum}}\dfrac{\pt}{\pt x_p}\left(\mathrm{m}^{pq}_{i}( y, z)\left( \dfrac{\pt u_{i}}{\pt x_q}+\dfrac{1}{\e}\dfrac{\pt u_{i}}{\pt y_q}+\dfrac{1}{\e \delta}\dfrac{\pt u_{i}}{\pt z_q}\right)\right)\right]\left( t, x,\dfrac{x}{\e},\dfrac{x}{\e \delta} \right) 
\\ & \ \ \ -\dfrac{1}{\e}\left[\overset{d}{\underset{p,q=1}{\sum}}\dfrac{\pt}{\pt y_p}\left(\mathrm{m}^{pq}_{i}( y, z)\left( \dfrac{\pt u_{i}}{\pt x_q}+\dfrac{1}{\e}\dfrac{\pt u_{i}}{\pt y_q}+\dfrac{1}{\e \delta}\dfrac{\pt u_{i}}{\pt z_q}\right)\right)\right]\left( t, x,\dfrac{x}{\e},\dfrac{x}{\e \delta} \right)
\\ & \ \ \ -\dfrac{1}{\e \delta}\left[\overset{d}{\underset{p,q=1}{\sum}}\dfrac{\pt}{\pt z_p}\left(\mathrm{m}^{pq}_{i}( y, z)\left( \dfrac{\pt u_{i}}{\pt x_q}+\dfrac{1}{\e}\dfrac{\pt u_{i}}{\pt y_q}+\dfrac{1}{\e \delta}\dfrac{\pt u_{i}}{\pt z_q}\right)\right)\right]\left( t, x,\dfrac{x}{\e},\dfrac{x}{\e \delta} \right) 
\\&=[(\e^{-2}\delta^{-2} \A_{zz}+\e^{-2}\delta^{-1} \A_{yz}+\e^{-1}\delta^{-1}\A_{xz}
\\& \ \ \ +\e^{-2} \A_{yy}+\e^{-1}\A_{xy}+\e^{0}\delta^{0} \A_{xx} )u_i]\left( t, x, \dfrac{x}{\e}, \dfrac{x}{\e \delta} \right)
\end{aligned}
\label{Aeue_intra}
\end{equation} 
with each operator defined by:
\begin{equation*}
 \begin{cases}\A_{s s}=-\overset{d}{\underset{p,q=1}{\sum}}\dfrac{\pt}{\pt s_p}\left( \mathrm{m}^{pq}_{i}( y,z)\dfrac{\pt}{\pt s_q} \right), \\ \A_{s h}=-\overset{d}{\underset{p,q=1}{\sum}}\dfrac{\pt}{\pt s_p}\left(\mathrm{m}^{pq}_{i}( y, z)\dfrac{\pt}{\pt h_q}\right)-\overset{d}{\underset{p,q=1}{\sum}}\dfrac{\pt}{\pt h_q}\left(\mathrm{m}^{pq}_{i}( y,z)\dfrac{\pt}{\pt s_p}\right) \quad \text{ if } s \neq h,
 \end{cases}
\end{equation*}
for  $s,h:=x,y,z.$ 

Now, we substitute the asymptotic expansion \eqref{ue_intra} of $u^{\e,\delta}_{i}$ into the operator developed \eqref{Aeue_intra} to obtain :
\begin{equation*}
\begin{aligned}
&\A_{\e,\delta} u_i^{\e,\delta}( t, x)  
\\&=[\e^{-2}\delta^{-2}\A_{zz} u_{i,0}+\e^{-2}\delta^{-1}\A_{yz} u_{i,0}+\e^{-2}\A_{yy} u_{i,0}+\e^{-1}\delta^{-2}\A_{zz} u_{i,1}+\delta^{-2} \A_{zz} u_{i,3}
\\ & \ \ +\e^{-1}\delta^{-1} \left(\A_{zz} u_{i,2}+\A_{yz} u_{i,1}+\A_{xz} u_{i,0} \right) +\delta^{-1}\left(\A_{zz} u_{i,4}+\A_{yz} u_{i,3}+\A_{xz} u_{i,1} \right)
\\ & \ \ +\e^{-1}\left(\A_{yz} u_{i,2}+\A_{yy} u_{i,1}+\A_{xz} u_{i,0} \right)
\\ & \ \ +\e^{0}\delta^{0}\left(\A_{zz} u_{i,5}+\A_{yz} u_{i,4}+\A_{yy} u_{i,3}+\A_{xz} u_{i,2}+\A_{xy} u_{i,1}+\A_{xx} u_{i,0} \right)]\left( t, x,\dfrac{x}{\e},\dfrac{x}{\e \delta} \right)+\cdots.
\end{aligned}
\end{equation*}
 Similarly, we have the boundary condition:
\begin{equation}
\begin{aligned}
\mathrm{M}_{i}^{\e,\delta} \nabla u^{\e,\delta}_{i} \cdot n &= \left[ \mathrm{M}_{i}^{\e,\delta}\nabla_x u_{i} +\e^{-1} \mathrm{M}_{i}^{\e,\delta}\nabla_y u_{i} + \e^{-1}\delta^{-1} \mathrm{M}_{i}^{\e,\delta}\nabla_z u_{i}\right]\cdot n, 
\end{aligned}
\label{Aeuemultibord}
\end{equation}
for $n:=n_{i},n_{z}.$ \\
Thus, we also substitute the asymptotic expansion \eqref{ue_intra} of $u^{\e,\delta}_{i}$ into the boundary condition equation \eqref{pbinintra} on $\Gamma^{y} $ and on $\Gamma^{z}:$
\begin{equation*}
\begin{aligned}
\mathrm{M}_{i}^{\e,\delta} \nabla u_i^{\e,\delta} \cdot n &=[\e^{0}\delta^{0} (\mathrm{M}_{i}\nabla_x u_{i,0}) \cdot n+\e (\mathrm{M}_{i}\nabla_x u_{i,1}) \cdot n+\e\delta (\mathrm{M}_{i}\nabla_x u_{i,2})\cdot n+\cdots]\left( t, x, \dfrac{x}{\e}, \dfrac{x}{\e \delta} \right) 
\\ & \ \ +[\e^{-1} (\mathrm{M}_{i}\nabla_y u_{i,0}) \cdot n+\e^{0}\delta^{0} (\mathrm{M}_{i}\nabla_y u_{i,1}) \cdot n+\delta (\mathrm{M}_{i}\nabla_y u_{i,2}) \cdot n
\\ & \ \ \ \ +\e (\mathrm{M}_{i}\nabla_y u_{i,3}) \cdot n+\e\delta (\mathrm{M}_{i}\nabla_y u_{i,4}) \cdot n + \cdots]\left( t, x, \dfrac{x}{\e}, \dfrac{x}{\e \delta} \right)
\\ & \ \ +[\e^{-1}\delta^{-1} (\mathrm{M}_{i}\nabla_z u_{i,0}) \cdot n+\delta^{-1} (\mathrm{M}_{i}\nabla_z u_{i,1}) \cdot n+\e^{0}\delta^{0} (\mathrm{M}_{i}\nabla_z u_{i,2}) \cdot n
\\ & \ \ \ \ +\e\delta^{-1} (\mathrm{M}_{i}\nabla_z u_{i,3}) \cdot n+\e (\mathrm{M}_{i}\nabla_z u_{i,4}) \cdot n+\e \delta (\mathrm{M}_{i}\nabla_z u_{i,5}) \cdot n+\cdots]\left( t, x, \dfrac{x}{\e}, \dfrac{x}{\e \delta} \right) 
\\&=[\e^{-1}\delta^{-1} (\mathrm{M}_{i}\nabla_z u_{i,0}) \cdot n+\e^{-1} (\mathrm{M}_{i}\nabla_y u_{i,0} )\cdot n+\delta^{-1} (\mathrm{M}_{i}\nabla_z u_{i,1}) \cdot n
\\& \ \ +\e^{0}\delta^{0}\left(\mathrm{M}_{i}\nabla_z u_{i,2}+\mathrm{M}_{i}\nabla_y u_{i,1}+\mathrm{M}_{i}\nabla_x u_{i,0} \right)\cdot n +\e\delta^{-1}\left(\mathrm{M}_{i}\nabla_z u_{i,3} \right)\cdot n
\\& \ \  +\e\left(\mathrm{M}_{i}\nabla_z u_{i,4}+\mathrm{M}_{i}\nabla_y u_{i,3}+\mathrm{M}_{i}\nabla_x u_{i,1} \right)\cdot n+\delta (\mathrm{M}_{i}\nabla_y u_{i,2}) \cdot n
\\& \ \  +\e \delta\left(\mathrm{M}_{i}\nabla_z u_{i,5}+\mathrm{M}_{i}\nabla_y u_{i,4}+\mathrm{M}_{i}\nabla_x u_{i,2} \right)\cdot n]\left( t, x, \dfrac{x}{\e}, \dfrac{x}{\e \delta} \right)+\cdots.
\end{aligned}
\end{equation*}
where $n$ represents the outward unit normal on $\Gamma^{y} $ or on $\Gamma^{z}$ $(n:=n_{i},n_{z}).$
Consequently, by equating the terms of the powers coefficients $\e^{\ell}\delta^m$ for the elliptic equations  and of the powers coefficients $\e^{\ell+1}\delta^{m+1}$ for the boundary conditions ($\ell,m=-2,-1,0$), we obtain the following systems: 
 \begin{equation}
 \begin{cases}
\A_{zz} u_{i,0}=0 \ \text{in} \ Z_{c}, \\ u_{i,0} \ z\text{-periodic}, \\ \mathrm{M}_{i}\nabla_z u_{i,0} \cdot n_{z}=0 \ \text{on} \ \Gamma^{z},  
 \end{cases}
 \label{Azzui0}
 \end{equation}
 
  \begin{equation}
 \begin{cases}
\A_{yy} u_{i,0}=0 \ \text{in} \ Y_{i}, \\ u_{i,0} \ y\text{-periodic}, \\ \mathrm{M}_{i}\nabla_y u_{i,0} \cdot n_{i}=0 \ \text{on} \ \Gamma^{y},  
 \end{cases}
 \label{Ayyui0}
 \end{equation}

\begin{equation}
 \begin{cases}
\A_{yz} u_{i,0}=0 \ \text{in} \ Z_{c}, \\ u_{i,0} \ y\text{- and }z\text{-periodic}, \\ \mathrm{M}_{i}\nabla_y u_{i,0} \cdot n_{i}=0 \ \text{on} \ \Gamma^{y}, \\ \mathrm{M}_{i}\nabla_z u_{i,0} \cdot n_{z}=0 \ \text{on} \ \Gamma^{z},  
 \end{cases}
 \label{Ayzui0}
 \end{equation}

 \begin{equation}
 \begin{cases}
\A_{zz} u_{i,1}=0 \ \text{in} \ Z_{c}, \\ u_{i,1} \ z\text{-periodic}, \\ \mathrm{M}_{i}\nabla_z u_{i,1} \cdot n_{z}=0 \ \text{on} \ \Gamma^{z},
 \end{cases}
 \label{Azzui1}
 \end{equation}

 \begin{equation}
 \begin{cases}
\A_{zz} u_{i,2}=-\A_{yz} u_{i,1}-\A_{xz} u_{i,0} \ \text{in} \ Z_{c},\\ u_{i,2} \ z\text{-periodic}, \\ \left(\mathrm{M}_{i}\nabla_z u_{i,2}+\mathrm{M}_{i}\nabla_y u_{i,1}+\mathrm{M}_{i}\nabla_x u_{i,0} \right)\cdot n_{z}=0 \ \text{on} \ \Gamma^{z}, 
 \end{cases}
 \label{Azzui2}
 \end{equation}
 
 \begin{equation}
 \begin{cases}
\A_{zz} u_{i,3}=0 \ \text{in} \ Z_{c}, \\ u_{i,3} \ z\text{-periodic}, \\ \left(\mathrm{M}_{i}\nabla_z u_{i,3} \right)\cdot n_{z}=0 \ \text{on} \ \Gamma^{z}, 
 \end{cases}
 \label{Azzui3}
 \end{equation}
 
 \begin{equation}
 \begin{cases}
\A_{zz} u_{i,4}=-\A_{yz} u_{i,3}-\A_{xz} u_{i,1} \ \text{in} \ Z_{c},\\ u_{i,4} \ y\text{- and }z\text{-periodic}, \\  \left(\mathrm{M}_{i}\nabla_z u_{i,4}+\mathrm{M}_{i}\nabla_y u_{i,3}+\mathrm{M}_{i}\nabla_x u_{i,1}\right)\cdot n_{i}=-\left( \pt_{t} v_0+\I_{ion}( v_0, {w}_0)-\I_{app}\right) \ \text{on} \ \Gamma^{y}, \\ \left(\mathrm{M}_{i}\nabla_z u_{i,4}+\mathrm{M}_{i}\nabla_y u_{i,3}+\mathrm{M}_{i}\nabla_x u_{i,1} \right)\cdot n_{z}=0 \ \text{on} \ \Gamma^{z}, 
 \end{cases}
 \label{Azzui4}
 \end{equation}

 \begin{equation}
 \begin{cases}
\A_{yz} u_{i,2}=-\A_{yy} u_{i,1}-\A_{xy} u_{i,0} \ \text{in} \ Z_{c},\\ u_{i,2} \ y\text{- and }z\text{-periodic}, \\  \left(\mathrm{M}_{i}\nabla_z u_{i,2}+\mathrm{M}_{i}\nabla_y u_{i,1}+\mathrm{M}_{i}\nabla_x u_{i,0} \right)\cdot n_{i}=0 \ \text{on} \ \Gamma^{y}, \\ \mathrm{M}_{i}\nabla_y u_{i,2} \cdot n_{z}=0 \ \text{on} \ \Gamma^{z},  
 \end{cases}
 \label{Ayzui2}
 \end{equation}
 
 \begin{equation}
 \begin{cases}
\A_{zz} u_{i,5}=-\A_{yz} u_{i,4}-\A_{yy} u_{i,3}-\A_{xz} u_{i,2}-\A_{xy} u_{i,1}-\A_{xx} u_{i,0} \ \text{in} \ Z_{c},\\ u_{i,5} \ z\text{-periodic}, \\ \left(\mathrm{M}_{i}\nabla_z u_{i,5}+\mathrm{M}_{i}\nabla_y u_{i,4}+\mathrm{M}_{i}\nabla_x u_{i,2} \right)\cdot n_{z}=0 \ \text{on} \ \Gamma^{z}.  
 \end{cases}
 \label{Azzui5}
 \end{equation}
 
 These systems \eqref{Azzui0}-\eqref{Azzui5} have a particular structure in the sense that their unknowns will be found iteratively.

 We will solve these nine problems \eqref{Azzui0}-\eqref{Azzui5} successively to determine the homogenized problem (based on the work \cite{doina} and \cite{ben}). The resolution is described as follows:
\begin{enumerate}
\item[\textbf{$\bullet$ Step 1}]We begin with the first problem \eqref{Azzui0} whose the following variational formulation: 
\begin{equation}
 \begin{cases}
\text{Find} \ \dot{u}_{i,0} \in \W_{per}(Z_{c}) \ \text{such that} \\ \dot{a}_{Z_{c}}(\dot{u}_{i,0},\dot{v})=\displaystyle \int_{\pt Z_{c}} \left( \mathrm{M}_i \nabla_z u_{i,0} \cdot n_{z}\right) v \ d\sigma_z,\ \forall \dot{v}\in \W_{per}(Z_{c}),
 \end{cases}
 \label{FvAzzui0}
\end{equation} 
with $\dot{a}_{Z_{c}}$ given by:
\begin{equation}
\dot{a}_{Z_{c}}(\dot{u},\dot{v})=\int_{Z_{c}} \mathrm{M}_i \nabla_z u \nabla_z  v \ dz, \ \forall u \in\dot{u}, \ \forall v \in\dot{v},\ \forall \dot{u},\ \forall \dot{v} \in \W_{per}(Z_{c})
\label{dotazi}
\end{equation}
and $$ \W_{per}(Z_{c})=H^1_{per}(Z_{c})/\R$$ is given by Definition \ref{w_per}. Similarly, we want to clarify the right hand side of the variational formulation \eqref{FvAzzui0}. By the definition of $\pt Z_{c}:= \pt_{\text{ext}} Z \cup \Gamma^{z},$ we use Proposition \ref{cond_Y_per} and the $z$-periodicity of $\mathrm{M}_i$ by taking account the boundary condition on $\Gamma^{z}$ to say that :
   \begin{align*}
 \displaystyle \int_{\pt Z_{c}} &\left(\mathrm{M}_i \nabla_z u_{i,0} \cdot n_{z}\right) v \ d\sigma_z 
 \\& \ \ \ =\int_{\pt_{\text{ext}} Z} \left(\mathrm{M}_i \nabla_z u_{i,0} \cdot n_{z}\right) v \ d\sigma_z+\int_{\Gamma^{z}} \left(\mathrm{M}_i \nabla_z u_{i,0} \cdot n_{z}\right) v \ d\sigma_z =0.
 \end{align*}
 Using Theorem \ref{cond_per_neu}, we obtain the existence and the uniqueness of solution $\dot{u}_{i,0}$ to the problem \eqref{FvAzzui0}. In addition, we have:
$$\norm {\dot{u}_{i,0}}_{\W_{per}(Z_{c})}=0.$$   So, $u_{i,0}$ is independent of the microscopic variable $z$. Thus, we deduce that: $$u_{i,0}(t,x,y,z)=u_{i,0}(t,x,y),\ \forall u_{i,0} \in\dot{u}_{i,0}.$$

\item[$\bullet$ \textbf{Step 2}] We now solve the second boundary value problem \eqref{Ayyui0} that is defined in $Y_{i}$. Its variational formulation is: 
\begin{equation}
 \begin{cases}
\text{Find} \ \dot{u}_{i,0} \in \W_{per}(Y_{i}) \ \text{such that} \\ \\ \dot{a}_{Y_{i}}(\dot{u}_{i,0},\dot{v})= \displaystyle \int_{\pt Y_i} \mathrm{M}_i \nabla_y u_{i,0} \cdot n_{i} \ v \ d\sigma_y \ \forall \dot{v}\in \W_{per}(Y_{i}),
 \end{cases}
 \label{FvAyyui0}
\end{equation} 
with $\dot{a}_ {Y_{i}}$ given by:
\begin{equation}
\dot{a}_ {Y_{i}}(\dot{u},\dot{v})=\int_{Y_{i}} \mathrm{M}_i \nabla_y u \nabla_y v dy, \ \forall u \in\dot{u}, \ \forall v \in\dot{v},\ \forall \dot{u},\ \forall \dot{v} \in \W_{per}(Y_{i})
\label{dotayi}
\end{equation}
and $ \W_{per}(Y_{i})$ given by Definition \ref{w_per}.
 
   Similarly, we want to clarify first the right hand side in the variational formulation \eqref{FvAyyui0}. By the definition of $\pt Y_i:=(\pt_{\text{ext}} Y \cap \pt Y_i)\cup \Gamma^{y},$ wwe use Proposition \ref{cond_Y_per} and the $y$-periodicity of $\mathrm{M}_i$ by taking account the boundary condition on $\Gamma^{y}$ to say that:  
   \begin{align*}
 \int_{\pt Y_{i}} &\mathrm{M}_{i} \nabla_y u_{i,0} \cdot n_{i}(y) v \ d\sigma_y 
 \\ & \ \ \ =\int_{\pt_{\text{ext}} Y  \cap \pt Y_i} \mathrm{M}_{i} \nabla_y u_{i,0} \cdot n_{i}(y) v \ d\sigma_y +\int_{\Gamma^{y}} \mathrm{M}_{i} \nabla_y u_{i,0} \cdot n_{i}(y) \ d\sigma_y=0.
 \end{align*}
 Therefore, we can apply Theorem \ref{cond_per_neu} to prove the existence and uniqueness of solution $\dot{u}_{i,0}.$ In addition, we have:
$$\norm {\dot{u}_{i,0}}_{\W_{per}(Y_{i})}=0.$$
Thus, we deduce that $u_{i,0}$ is also independent of the mesoscopic variable $y.$ Consequently, the third boundary value problem \eqref{Ayzui0} is satisfied automatically.\\
Next, we solve the fourth problem \eqref{Azzui1} by the same process of the first step. So, we deduce that $u_{i,1}$ is independent of $z$. Finally, we have: 
$$u_{i,0}( t, x, y, z)=u_{i,0}(t, x) \text{ and } u_{i,1}(t,x,y,z)=u_{i,1}(t,x,y).$$
\begin{rem} Since $u_{i,0}$ is independent of $y$ and $z$ then it does not oscillate "rapidly". This is why now expect $u_{i,0}$ to be the "homogenized solution". To find the homogenized equation, it is sufficient to find an equation in $\Omega$ satisfied by $u_{i,0}$ independent on $y$ and $z.$
\end{rem}

\item[$\bullet$ \textbf{Step 3}] We solve the fifth problem \eqref{Azzui2}. Taking into account the form of $u_{i,0}$ and $u_{i,1},$ system \eqref{Azzui2} can be rewritten as:
\begin{equation}
 \begin{cases}
\A_{zz} u_{i,2}=\overset{d}{\underset{p,q=1}{\sum}}\dfrac{\pt \mathrm{m}^{pq}_{i}}{\pt z_p}\left(\dfrac{\pt u_{i,1}}{\pt y_q}+\dfrac{\pt u_{i,0}}{\pt x_q}\right) \ \text{in} \ Z_{c},\\ u_{i,2} \ z\text{-periodic}, \\ \left(\mathrm{M}_{i}\nabla_z u_{i,2}+\mathrm{M}_{i}\nabla_y u_{i,1}+\mathrm{M}_{i}\nabla_x u_{i,0} \right)\cdot n_{z}=0 \ \text{on} \ \Gamma^{z}, 
 \end{cases}
 \label{Azzui2new}
 \end{equation}
 Its variational formulation is:
 \begin{equation}
 \begin{cases}
\text{Find} \ \dot{u}_{i,2} \in \W_{per}(Z_{c}) \ \text{such that} \\ \dot{a}_{Z_{c}}(\dot{u}_{i,2},\dot{v})=(F_2,\dot{v})_{(\W_{per}(Z_{c}))',\W_{per}(Z_{c})} \ \forall \dot{v}\in \W_{per}(Z_{c}),
 \end{cases}
 \label{FvAzzui2}
\end{equation}
 with $\dot{a}_{Z_{c}}$ given by \eqref{dotazi} and $F_2$ defined by:
\begin{equation}
(F_2,\dot{v})_{(\W_{per}(Z_{c}))',\W_{per}(Z_{c})}=-\overset{d}{\underset{p,q=1}{\sum}}\left(\dfrac{\pt u_{i,1}}{\pt y_q}+\dfrac{\pt u_{i,0}}{\pt x_q}\right) \int_{Z_{c}} \mathrm{m}^{pq}_{i}( t, y, z)\dfrac{\pt v}{\pt z_p}dz,
\label{F2i}
\end{equation}
for all $v \in\dot{v}$ and $\dot{v} \in \W_{per}(Z_{c}).$
  
  Note that $F_2$ belongs to $(\W_{per}(Z_{c}))'.$ Then, Theorem \ref{cond_per_neu} gives a unique solution $\dot{u}_{i,2} \in \W_{per}(Z_{c})$ of the problem \eqref{Azzui2new}-\eqref{F2i}. \\
 Thus, the linearity of terms in the right of equation \eqref{Azzui2new} suggests to look for $\dot{u}_{i,2}$ under the following form:
\begin{equation}
\dot{u}_{i,2}=\dot{\theta}_{i}(z)\cdot \left(\nabla_y \dot{u}_{i,1}+\nabla_x \dot{u}_{i,0}\right) \ \text{in} \ \W_{per}(Z_{c}),
\label{dotui2}
\end{equation}
with the corrector function $\dot{\theta}_{i}^q$ (i.e the components of the function  $\dot{\theta}_{i}$) satisfies the $\delta$-cell problem:
\begin{equation}
\begin{cases}
\A_{zz} \dot{\theta}_{i}^q =\overset{d}{\underset{p=1}{\sum}}\dfrac{\pt \mathrm{m}^{pq}_{i}}{\pt z_p}(y,z) \ \text{in} \ Z_{c},\\ \dot{\theta}_{i}^q \ y\text{- and }z\text{-periodic}, \\ \mathrm{M}_{i}\nabla_z \dot{\theta}_{i}^q \cdot n_{z}=-(\mathrm{M}_{i}e_q)\cdot n_{z} \ \text{on} \ \Gamma^{z}, 
 \end{cases}
 \label{Azzthetai}
 \end{equation}
 for $e_q, \, q=1,\dots,d,$ the standard canonical basis in $\R^d.$ Moreover, we can choose a representative element $\theta_{i}^q$ of the class $\dot{\theta}_{i}^q$ which satisfy the following variational formulation:
 \begin{equation}
 \begin{cases}
\text{Find} \ \theta_{i}^q \in W_{\#}(Z_{c}) \ \text{such that} \\ \\ a_{Z_{c}} (\theta_{i}^q,v) =-\overset{d}{\underset{p=1}{\sum}}\displaystyle\int_{Z_{c}} \mathrm{m}^{pq}_{i}( t, y, z) \dfrac{\pt v}{\pt z_p} \ dz, \ \forall v \in W_{\#}(Z_{c}),
 \end{cases}
 \label{FvAzzthetai}
\end{equation}
with $W_{\#}(Z_{c})$ given by the expression \eqref{W}. The condition of Theorem \ref{cond_per_neu} is imposed to guarantee the existence and uniqueness of the solution of the problem \eqref{Azzthetai}-\eqref{FvAzzthetai}. 
Thus, by the form $\dot{u}_{i,2}$ given by the expression \eqref{dotui2}, the solution $u_{i,2}$  can be represented by the following ansatz:
\begin{equation} 
u_{i,2}(t,x,y,z)=\theta_{i}(z)\cdot \left(\nabla_y u_{i,1}(t,x,y)+\nabla_x u_{i,0}(t,x)\right)+\tilde{u}_{i,2}(t,x,y) \ \text{with} \ u_{i,2} \in \dot{u}_{i,2},
\label{ui2}
\end{equation}
and $\tilde{u}_{i,2}$ is a constant with respect to $z$ (i.e $\tilde{u}_{i,2} \in \dot{0}$ in $\W_{\#}(Z_{c})$).

 Next, we pass to the sixth problem \eqref{Azzui3} by the same strategy of the first step. We obtain that $u_{i,3}$ is independent of $z$ and we have:
 \begin{equation*}
 u_{i,3}(t,x,y,z)=u_{i,3}(t,x,y).
 \label{ui3}
 \end{equation*}

\item[$\bullet$ \textbf{Step 4}] We now solve the seventh boundary value problem \eqref{Azzui4}. Taking into account the form of $u_{i,3}$ and $u_{i,1},$ we can rewrite this problem as follows: 
\begin{equation}
 \begin{cases}
\A_{zz} u_{i,4}=\overset{d}{\underset{p,q=1}{\sum}}\dfrac{\pt \mathrm{m}^{pq}_{i}}{\pt z_p}\left(\dfrac{\pt u_{i,3}}{\pt y_q}+\dfrac{\pt u_{i,1}}{\pt x_q}\right) \ \text{in} \ Z_{c},\\ u_{i,4} \ y\text{- and }z\text{-periodic}, \\  \left(\mathrm{M}_{i}\nabla_z u_{i,4}+\mathrm{M}_{i}\nabla_y u_{i,3}+\mathrm{M}_{i}\nabla_x u_{i,1} \right)\cdot n_{z}=0 \ \text{on} \ \Gamma^{z}. 
 \end{cases}
 \label{Azzui4new}
 \end{equation}

 Its variational formulation  is:
\begin{equation}
 \begin{cases}
\text{Find} \ \dot{u}_{i,4} \in \W_{per}(Z_{c}) \ \text{such that} \\ \dot{a}_{Z_{c}}(\dot{u}_{i,4},\dot{v})=(F_4,\dot{v})_{(\W_{per}(Z_{c}))',\W_{per}(Z_{c})} \ \forall \dot{v}\in \W_{per}(Z_{c}),
 \end{cases}
 \label{FvAzzui4}
\end{equation}
with $\dot{a}_{Z_{c}}$ given by \eqref{dotazi} and $F_4$ defined by:
\begin{equation}
(F_4,\dot{v})_{(\W_{per}(Z_{c}))',\W_{per}(Z_{c})}=-\overset{d}{\underset{p,q=1}{\sum}}\left(\dfrac{\pt u_{i,3}}{\pt y_q}+\dfrac{\pt u_{i,1}}{\pt x_q}\right) \int_{Z_{c}} \mathrm{m}^{pq}_{i}( t, y, z)\dfrac{\pt v}{\pt z_p}dz,
\label{F4i}
\end{equation}
for all $v \in\dot{v}$ and $\dot{v} \in \W_{per}(Z_{c}).$\\
The problem \eqref{Azzui4new}-\eqref{F4i} is well-posed according to Theorem \ref{cond_per_neu} under the compatibility condition: 
\begin{equation*}
(F_4,1)_{(\W_{per}(Z_{c}))',\W_{per}(Z_{c})}=0.
\end{equation*}
 This implies that problem \eqref{Azzui4} has a unique periodic solution up to a constant. Thus, the linearity of terms in the right hand side of equation \eqref{Azzui4new} suggests to look for $u_{i,4}$ under the following form:
\begin{equation} 
u_{i,4}(t,x,y,z)=\theta_{i}(z)\cdot \left(\nabla_y u_{i,3}(t,x,y)+\nabla_x u_{i,1}(x)\right)+\tilde{u}_{i,4}(t,x,y) \ \text{with} \ u_{i,4} \in \dot{u}_{i,4},
\label{ui4}
\end{equation}
where $\tilde{u}_{i,4}$ is a constant with respect to $z$ and $\theta_{i}$ satisfies problem \eqref{Azzthetai}.
\item[$\bullet$ \textbf{Step 5}] We consider the eighth problem \eqref{Ayzui2}:
\begin{equation*}
 \begin{cases}
\A_{yz} u_{i,2}=-\A_{yy} u_{i,1}-\A_{xy} u_{i,0} \ \text{in} \ Z_{c},\\ u_{i,2} \ z\text{-periodic}, \\  \left(\mathrm{M}_{i}\nabla_z u_{i,2}+\mathrm{M}_{i}\nabla_y u_{i,1}+\mathrm{M}_{i}\nabla_x u_{i,0} \right)\cdot n_{i}=0 \ \text{on} \ \Gamma^{y}, \\ \mathrm{M}_{i}\nabla_y u_{i,2} \cdot n_{z}=0 \ \text{on} \ \Gamma^{z}.  
 \end{cases}
 \end{equation*} 
  
 Taking into account the form of $u_{i,0}$ and $u_{i,1},$ we can rewrite the first equation as follows: 
\begin{equation*}
\A_{yz} u_{i,2} =\overset{d}{\underset{p,q=1}{\sum}}\dfrac{\pt}{\pt y_p}\left(\mathrm{m}^{pq}_{i}( y, z)\dfrac{\pt u_{i,1}}{\pt y_q}\right) +\overset{d}{\underset{p,q=1}{\sum}}\dfrac{\pt \mathrm{m}^{pq}_{i}}{\pt y_p}( y,z)\dfrac{\pt u_{i,0}}{\pt x_q}. 
\label{Ayzui2new}
\end{equation*}
To find the explicit form of $u_{i,1}$, we will follow the following steps:  
  First, we integrate over $Z_{c}$ the above equation as follows:  
\begin{equation}
\begin{aligned}
 &-\overset{d}{\underset{p,q=1}{\sum}}\int_{Z_{c}} \dfrac{\pt}{\pt y_p}\left(\mathrm{m}^{pq}_{i}( y, z)\dfrac{\pt u_{i,2}}{\pt z_q}\right)\ dz -\overset{d}{\underset{p,q=1}{\sum}}\int_{Z_{c}} \dfrac{\pt}{\pt z_p}\left(\mathrm{m}^{pq}_{i}( y, z)\dfrac{\pt u_{i,2}}{\pt y_q}\right) \ dz
 \\ & \qquad \qquad  =\overset{d}{\underset{p,q=1}{\sum}}\int_{Z_{c}}\dfrac{\pt}{\pt y_p}\left(\mathrm{m}^{pq}_{i}( y, z)\dfrac{\pt u_{i,1}}{\pt y_q}\right) +\overset{d}{\underset{p,q=1}{\sum}}\int_{Z_{c}}\dfrac{\pt \mathrm{m}^{pq}_{i}}{\pt y_p}( y, z)\dfrac{\pt u_{i,0}}{\pt x_q} \ dz.
 \end{aligned}
 \label{Int_ZAyzui2new}
 \end{equation}
 We denote by $E_i$ with $i = 1,\dots, 4$ the terms of the previous equation which is rewritten as follows (to respect the order):
\begin{equation*}
E_1+E_2=E_3+E_4.
\end{equation*}
Next, we use the divergence formula for the second term $E_2$ together with Proposition \ref{cond_Y_per} and the boundary condition on $\Gamma^{z}$  to obtain:
\begin{align*}
E_2 & =-\int_{\pt Z_{c}} \mathrm{M}_{i}\nabla_y u_{i,2} \cdot n_z \ d\sigma_z
\\& =-\int_{\pt_{\text{ext}} Z} \mathrm{M}_{i}\nabla_y u_{i,2} \cdot n_z \ d\sigma_z-\int_{\Gamma^{z}} \mathrm{M}_{i}\nabla_y u_{i,2} \cdot n_z \ d\sigma_z=0.
 \end{align*}
Now, we replace $u_{i,2}$ by its expression \eqref{ui2} in the first term $E_1$ to obtain the following:
\begin{equation*}
 E_1=-\overset{d}{\underset{p,q=1}{\sum}}\int_{Z_{c}} \dfrac{\pt}{\pt y_p}\left(\mathrm{m}^{pq}_{i}( y, z)\left( \overset{d}{\underset{k=1}{\sum}}\dfrac{\pt \theta^k_{i}}{\pt z_q}\left(\dfrac{\pt u_{i,1}}{\pt y_k}+\dfrac{\pt u_{i,0}}{\pt x_k}\right)\right)\right)\ dz.
\end{equation*}
By permuting the index in the right hand side of the equation \eqref{Int_ZAyzui2new}, we obtain:
\begin{align*}
& E_3 =\overset{d}{\underset{p,k=1}{\sum}}\int_{Z_{c}}\dfrac{\pt}{\pt y_p}\left(\mathrm{m}^{pk}_{i}( y, z)\dfrac{\pt u_{i,1}}{\pt y_k}\right),
\\& E_4=\overset{d}{\underset{p,k=1}{\sum}}\int_{Z_{c}}\dfrac{\pt \mathrm{m}^{pk}_{i}}{\pt y_p}( y, z)\dfrac{\pt u_{i,0}}{\pt x_k} \ dz.
\end{align*}
Finally, we obtain an equation for the mesoscopic scale (independent of $z$) satisfied by $u_{i,1}:$
\begin{align*}
 &- \overset{d}{\underset{p,k=1}{\sum}}\dfrac{\pt}{\pt y_p}\left(\dfrac{1}{\abs{Z}}\overset{d}{\underset{q=1}{\sum}}\left[\int_{Z_{c}} \left(\mathrm{m}^{pk}_{i}+ \mathrm{m}^{pq}_{i}\dfrac{\pt \theta_{i}^{k}}{\pt z_q}\right) \ dz \right]\dfrac{\pt u_{i,1}}{\pt y_k} \right)
 \\& \qquad =\overset{d}{\underset{p,k=1}{\sum}}\dfrac{\pt}{\pt y_p}\left(\dfrac{1}{\abs{Z}}\overset{d}{\underset{q=1}{\sum}}\left[\int_{Z_{c}} \left(\mathrm{m}^{pk}_{i}+ \mathrm{m}^{pq}_{i}\dfrac{\pt \theta_{i}^{k}}{\pt z_q}\right) \ dz \right]\right)\dfrac{\pt u_{i,0}}{\pt x_k}. 
\end{align*} 

 Similarly, we replace $u_{i,2} $ by its form \eqref{ui2} in the boundary condition on $\Gamma^{y}$ then we integrate over $Z_{c}$ to obtain another condition satisfied by $u_{i,1}.$  Then, we obtain a mesoscopic problem defined on the unit cell portion $Y_i$ and satisfied by $u_{i,1}$ as follows:
\begin{equation}
 \begin{cases}
\ \B_{yy} u_{i,1}=\overset{d}{\underset{p,k=1}{\sum}} \dfrac{\pt \widetilde{\mathbf{m}}^{pk}_{i} }{\pt y_p} \dfrac{\pt u_{i,0}}{\pt x_k} \ \text{in} \ Y_{i},\\ \\  \ \left( \widetilde{\mathbf{M}}_i \nabla_y u_{i,1} + \widetilde{\mathbf{M}}_i \nabla_x u_{i,0}\right)\cdot n_{i}=0 \ \text{ on } \ \Gamma^{y},  
 \end{cases}
 \label{Byyui1}
 \end{equation}
with the operator $\B_{yy}$ (homogenized operator with respect to $z$) defined by:
 \begin{equation}
\B_{yy}=-\overset{d}{\underset{p,k=1}{\sum}}\dfrac{\pt}{\pt y_p}\left(\widetilde{\mathbf{m}}^{pk}_{i}(y)\dfrac{\pt}{\pt y_k}\right),
\label{Byyhome}
\end{equation}
 
 where with the coefficients of the (homogenized with respect to z) conductivity matrices $\widetilde{\mathbf{M}}_i=(\widetilde{\mathbf{m}}^{pk}_i)_{1\leq p,k \leq d}$ defined by:   
\begin{equation}
 \widetilde{\mathbf{m}}^{pk}_i(y)=\dfrac{1}{\abs{Z}}\overset{d}{\underset{q=1}{\sum}}\displaystyle \int_{Z_{c}} \left(\mathrm{m}^{pk}_{i}+ \mathrm{m}^{pq}_{i}\dfrac{\pt \theta_{i}^{k}}{\pt z_q}\right) \ dz, \ \forall p,k=1,\dots,d.
 \label{Mtt_i}
\end{equation}
Note that the $y$-periodicity of function $\widetilde{\mathbf{m}}^{pk}_i$ comes from the fact that the coefficients of conductivity matrix $\mathrm{M}_{i}$ and of the function $\theta_{i}$ are $y$-periodic.
\begin{rem}The operator $\B_{yy}$ has the same properties of the \textbf{homogenized} operator \eqref{Bxxhome} for the extracellular problem. At this point, we deduce that this method is used to homogenize the problem with respect to $z$ and then with respect to $y$. We remark also that  allows to obtain the effective properties at $\delta$-structural level and which become the input values in order to find the effective behavior of the cardiac tissue.
\end{rem}  
   
  Now, we prove the existence and uniqueness of solution of the problem \eqref{Byyui1} defined in $Y_i$. Consider the variational formulation
  of problem \eqref{Byyui1}:
 \begin{equation}
 \begin{cases}
\text{Find} \ \dot{u}_{i,1} \in \W_{per}(Y_{i}) \ \text{such that} \\ \dot{b}_{Y_{i}}(\dot{u}_{i,1},\dot{v})=(F_1,\dot{v})_{(\W_{per}(Y_{i}))',\W_{per}(Y_{i})} \ \forall \dot{v}\in \W_{per}(Y_{i}),
 \end{cases}
 \label{FvByyui1}
\end{equation}
 with $\dot{b}_Y$ given by:
 \begin{equation}
\dot{b}_ {Y_{i}}(\dot{u},\dot{v})=\int_{Y_{i}} \widetilde{\mathbf{M}}_i \nabla_y u \nabla_y v dy, \ \forall u \in\dot{u}, \ \forall v \in\dot{v},\ \forall \dot{u},\ \forall \dot{v} \in \W_{per}(Y_{i})
\label{dotbyi}
\end{equation}
and $F_1$ defined by:
\begin{equation}
(F_1,\dot{v})_{(\W_{per}(Y_{i}))',\W_{per}(Y_{i})}=-\overset{d}{\underset{p,k=1}{\sum}}\dfrac{\pt u_{i,0}}{\pt x_k} \int_{Y_{i}} \widetilde{\mathbf{m}}_{i}^{pk}(y)\dfrac{\pt v}{\pt y_p}dy, \ \forall v \in\dot{v}, \ \forall \dot{v} \in \W_{per}(Y_{i}).
\label{F1i}
\end{equation}
 The linear form $F_1$ belongs to $(\W_{per}(Y_i))'.$ Thus, there exists a unique solution $\dot{u}_{i,1} \in \W_{per}(Y_{i})$ of problem \eqref{FvByyui1}-\eqref{F1i}. \\
 Finally, the linearity of terms in the right of the equation \eqref{Byyui1} suggests to look for $\dot{u}_{i,2}$ under the following form:
 \begin{equation}
\dot{u}_{i,1}=\dot{\chi}(y)\cdot\nabla_x\dot{u}_{i,0} \ \text{in} \ \W_{per}(Y_{i}),
\label{dotui1}
\end{equation}
with each element of the corrector function $\dot{\chi}_{i}=\left(\dot{\chi}_{i}^k\right)_{k=1,\dots,d}$ satisfies the following $\e$-cell problem:
\begin{equation}
\begin{cases}
\B_{yy} \dot{\chi}_{i}^k =\overset{d}{\underset{p=1}{\sum}}\dfrac{\pt \widetilde{\mathbf{m}}_{i}^{pk}}{\pt y_p} \ \text{in} \ Y_{i},\\ \\   \widetilde{\mathbf{M}}_{i} \nabla_y \dot{\chi}_{i}^k \cdot n_{i}=-\left(\widetilde{\mathbf{M}}_{i} e_k \right)\cdot n_{i} \ \text{on} \ \Gamma^{y},  
 \end{cases}
 \label{Byychii}
 \end{equation}
 for $e_k,k=1,\dots,d,$ the standard canonical basis in $\R^d.$ Moreover, we can choose a representative element $\chi_{i}^k$ of the class $\dot{\chi}_{i}^k$ which satisfy the following variational formulation:
 \begin{equation}
 \begin{cases}
\text{Find} \ \chi_{i}^k \in W_{\#}(Y_{i}) \ \text{such that} \\ \\ \dot{b}_ {Y_{i}}(\chi_{i}^k,\dot{v}) =-\displaystyle \overset{d}{\underset{p=1}{\sum}} \int_{Y_{i}} \widetilde{\mathbf{m}}_{i}^{pk}(y)  \dfrac{\pt w}{\pt y_p} dy,  \ \forall w \in W_{\#}(Y_{i}),
 \end{cases}
 \label{Fvchii}
\end{equation}
with $\dot{b}_{Y_{i}}$ given by \eqref{dotbyi}. Thus, we prove the existence and uniqueness of the solution $\chi_{i}^k$ of the problem \eqref{Byychii} using Theorem \ref{cond_per_neu}.\\
So, by the form of $\dot{u}_{i,1}$  given by \eqref{dotui1}, the solution $u_{i,1}$ of the problem \eqref{Byyui1} can be represented by the following ansatz:
\begin{equation} 
u_{i,1}(t,x,y)=\chi_{i}(y)\cdot \nabla_x u_{i,0}(t,x)+\tilde{u}_{i,1}(t,x) \ \text{avec} \ u_{i,1} \in \dot{u}_{i,1},
\label{ui1}
\end{equation}
where $\tilde{u}_{i,1}$ is a constant with respect to $y,$ (i.e $\tilde{u}_{i,1} \in \dot{0}$ in $\W_{per}(Y_{i})$).

\item[ $\bullet$ \textbf{Last step}] Our interest is the last boundary value problem \eqref{Azzui5}. We have 
\begin{equation*}
\begin{aligned}
-\A_{yz} &u_{i,4}-\A_{yy} u_{i,3}-\A_{xz} u_{i,2}-\A_{xy} u_{i,1}-\A_{xx} u_{i,0} 
\\ &=\overset{d}{\underset{p,q=1}{\sum}}\dfrac{\pt}{\pt y_p}\left(\mathrm{m}^{pq}_{i}( y, z)\left(\dfrac{\pt u_{i,4}}{\pt z_q}+\dfrac{\pt u_{i,3}}{\pt y_q}+\dfrac{\pt u_{i,1}}{\pt x_q}\right)\right) 
\\& \ \ +\overset{d}{\underset{p,q=1}{\sum}}\dfrac{\pt}{\pt z_p}\left(\mathrm{m}^{pq}_{i}( y, z)\left(\dfrac{\pt u_{i,4}}{\pt y_q}+\dfrac{\pt u_{i,2}}{\pt x_q}\right)\right)
\\ & \ \ +\overset{d}{\underset{p,q=1}{\sum}}\dfrac{\pt}{\pt x_p}\left(\mathrm{m}^{pq}_{i}( y, z)\left(\dfrac{\pt u_{i,2}}{\pt z_q}+\dfrac{\pt u_{i,1}}{\pt y_q}+\dfrac{\pt u_{i,0}}{\pt x_q}\right)\right). 
\end{aligned}
\end{equation*}
Note that, the variational formulation of system \eqref{Azzui5} can be written as follows:
\begin{equation}
 \begin{cases}
\text{Find} \ \dot{u}_{i,5} \in \W_{per}(Z_{c}) \ \text{such that} \\ \dot{a}_{Z_{c}}(\dot{u}_{i,5},\dot{v})=(F_5,\dot{v})_{(\W_{per}(Z_{c}))',\W_{per}(Z_{c})} \ \forall \dot{v}\in \W_{per}(Z_{c}),
 \end{cases}
 \label{FvAzzui5}
\end{equation}
with $\dot{a}_{Z_{c}}$ given by \eqref{dotazi} and $F_5$ defined by
\begin{equation}
\begin{aligned}
&(F_5,\dot{v})_{(\W_{per}(Z_{c}))',\W_{per}(Z_{c})} 
\\& \ \ =\int_{\Gamma^{z}}\left[ \left(\mathrm{M}_{i}\nabla_z u_{i,5}+\mathrm{M}_{i}\nabla_y u_{i,4}+\mathrm{M}_{i}\nabla_x u_{i,2} \right)\cdot n_{z}\right]v \ d\sigma_z 
\\& \ \ +\overset{d}{\underset{p,q=1}{\sum}} \int_{Z_{c}} \dfrac{\pt}{\pt y_p}\left(\mathrm{m}^{pq}_{i}( y, z)\left(\dfrac{\pt u_{i,4}}{\pt z_q}+\dfrac{\pt u_{i,3}}{\pt y_q}+\dfrac{\pt u_{i,0}}{\pt x_q}\right) \right) v dz
\\& \ \ -\overset{d}{\underset{p,q=1}{\sum}}\int_{Z_{c}} \mathrm{m}^{pq}_{i}( y, z)\left( \dfrac{\pt u_{i,4}}{\pt y_q}+\dfrac{\pt u_{i,2}}{\pt x_q}\right)\dfrac{\pt v}{\pt z_p} dz
\\ & \ \ +\overset{d}{\underset{p,q=1}{\sum}} \int_{Z_{c}} \dfrac{\pt}{\pt x_p}\left(\mathrm{m}^{pq}_{i}( y, z)\left(\dfrac{\pt u_{i,2}}{\pt z_q}+\dfrac{\pt u_{i,1}}{\pt y_q}+\dfrac{\pt u_{i,0}}{\pt x_q}\right) \right) v dz, \ \  \forall v \in\dot{v}, \ \forall \dot{v} \in \W_{per}(Z_{c}).
\end{aligned}
\label{F5i}
\end{equation}
 The aim is to find the homogenized equation in $\Omega.$ Firstly, we will homogenize the problem \eqref{Azzui5} with respect to $z.$ Next, we homogenize the last one with respect to $y$ using the explicit forms of previous solutions. Finally, we obtain the corresponding homogenized model.

 Firstly, the problem \eqref{FvAzzui5}-\eqref{F5i} defined in $Z_{c}$ is well-posed if and only if $F_5$ belongs to $(\W_{per}(Z_{c}))',$ i.e,
\begin{equation*}
(F_5,1)_{(\W_{per}(Z_{c}))',\W_{per}(Z_{c})}=0 
\end{equation*}
which equivalent to:
\begin{equation*}
\begin{aligned}
&  -\dfrac{1}{\abs{Z}}\overset{d}{\underset{p,q=1}{\sum}} \int_{Z_{c}} \dfrac{\pt}{\pt y_p}\left(\mathrm{m}^{pq}_{i}( y, z)\left(\dfrac{\pt u_{i,4}}{\pt z_q}+\dfrac{\pt u_{i,3}}{\pt y_q}+\dfrac{\pt u_{i,1}}{\pt x_q}\right) \right)dz
\\ & \ \ =\dfrac{1}{\abs{Z}}\overset{d}{\underset{p,q=1}{\sum}} \int_{Z_{c}} \dfrac{\pt}{\pt x_p}\left(\mathrm{m}^{pq}_{i}( y, z)\left(\dfrac{\pt u_{i,2}}{\pt z_q}+\dfrac{\pt u_{i,1}}{\pt y_q}+\dfrac{\pt u_{i,0}}{\pt x_q}\right) \right)dz.
\end{aligned}
\end{equation*}
  In addition, we replace $u_{i,4}$ by its expression \eqref{ui4} into the above condition and into the boundary condition equation on $\Gamma^{y}$ satisfied by $u_{i,4}.$ Then, we obtain that $u_{i,3}$ satisfies the following problem  defined in $Y_{i}$
 \begin{equation}
 \begin{cases}
\  \B_{yy} u_{i,3} =-\B_{xy} u_{i,1} -\B_{xx}u_{i,0} \ \text{in} \ Y_{i},\\ \\ \ \left(\widetilde{\mathbf{M}}_i \nabla_y u_{i,3}  +\widetilde{\mathbf{M}}_i\nabla_x u_{i,1}\right)\cdot n_{i}=-\left( \pt_{t} v_0+\I_{ion}(v_0,{w}_0)-\I_{app}\right) \text{ on } \Gamma^{y},  
 \end{cases}
 \label{Byyui3}
 \end{equation}
 with  $\B_{xy}:=-\nabla_x \cdot \left( \widetilde{\mathbf{M}}_i \nabla_y \right)-\nabla_y \cdot \left( \widetilde{\mathbf{M}}_i \nabla_x \right)$.\\
Consequently, system \eqref{Byyui3} have the following variational formulation:
\begin{equation}
 \begin{cases}
\text{Find} \ u_{i,3} \in \W_{per}(Y_{i}) \ \text{such that} \\ \dot{b}_{Y_{i}}(u_{i,3},\dot{w})=(F_3,\dot{w})_{(\W_{per}(Y_{i}))',\W_{per}(Y_{i})} \ \forall \dot{w}\in \W_{per}(Y_{i}),
 \end{cases}
 \label{FvByyui3}
\end{equation}
with $\dot{b}_{Y_i}$ given by \eqref{dotbyi} and $F_3$ defined by:
\begin{equation}
\begin{aligned}
&(F_3,\dot{w})_{(\W_{per}(Y_{i}))',\W_{per}(Y_{i})} 
\\ &=\displaystyle \int_{\Gamma^{y}}\left(\widetilde{\mathbf{M}}_i \nabla_y u_{i,3}  +\widetilde{\mathbf{M}}_i\nabla_x u_{i,1}\right)\cdot n_{i}w \ d\sigma_y -\overset{d}{\underset{p,k=1}{\sum}}  \displaystyle \int_{Y_{i}} \widetilde{\mathbf{m}}_i^{pk}\dfrac{\pt u_{i,1}}{\pt x_k}\dfrac{\pt w}{\pt y_p} dy 
\\& +\overset{d}{\underset{p,k=1}{\sum}} \displaystyle \int_{Y_{i}}\dfrac{\pt}{\pt x_p}\left(\widetilde{\mathbf{m}}_i^{pk}\left(\dfrac{\pt u_{i,1}}{\pt y_k}+\dfrac{\pt u_{i,0}}{\pt x_k}\right) \right) wdy, 
\end{aligned}
\label{F3i}
\end{equation}
for all $w \in\dot{w}, \ \dot{w} \in \W_{per}(Y_{i}).$
 
 Observe that problem \eqref{Byyui3}-\eqref{F3i} is well-posed if and only if $F_3$ belongs to $(\W_{per}Y))',$ which means
 \begin{equation*}
 (F_3,1)_{(\W_{per}(Y_{i}))',\W_{per}(Y_{i})}= 0
\end{equation*}  
which gives:
\begin{align*}
 -\overset{d}{\underset{p,k=1}{\sum}} \displaystyle \int_{Y_{i}}\dfrac{\pt}{\pt x_p}\left(\widetilde{\mathbf{m}}_i^{pk}\left(\dfrac{\pt u_{i,1}}{\pt y_k}+\dfrac{\pt u_{i,0}}{\pt x_k}\right) \right) dy= -\abs{\Gamma^{y}} \left( \pt_{t} v_0+\I_{ion}(v_0,{w}_0)-\I_{app}\right).
\end{align*}

 Next, we replace $u_{i,1}$ by its form \eqref{ui1} in the above condition. Then, we obtain:
 \begin{align*}
 &-\overset{d}{\underset{p,k=1}{\sum}}  \displaystyle \int_{Y_{i}} \dfrac{\pt}{\pt x_p}\left(\widetilde{\mathbf{m}}_i^{pk} \left(\overset{d}{\underset{q=1}{\sum}}\dfrac{\pt \chi_{i}^q}{\pt y_k}(y)\dfrac{\pt u_{i,0}}{\pt x_q}+\dfrac{\pt u_{i,0}}{\pt x_k}\right)\right) dy
 \\& \qquad\qquad = -\abs{\Gamma^{y}} \left( \pt_{t} v_0+\I_{ion}(v_0,{w}_0)-\I_{app}\right)
\end{align*}
 
 By expanding the sum and permuting the index, we obtain
\begin{align*}
& -\overset{d}{\underset{p,q=1}{\sum}}  \displaystyle \int_{Y_{i}} \dfrac{\pt}{\pt x_p}\left[\left( \overset{d}{\underset{k=1}{\sum}}\widetilde{\mathbf{m}}_i^{pk} \dfrac{\pt \chi_{i}^q}{\pt y_k}(y)+\widetilde{\mathbf{m}}_i^{pq}\right)\dfrac{\pt u_{i,0}}{\pt x_q}\right] dy
\\& \qquad\qquad = -\abs{\Gamma^{y}} \left( \pt_{t} v_0+\I_{ion}(v_0,{w}_0)-\I_{app}\right).
\end{align*}
Then, the function $u_{i,0}$ satisfies the following problem:
\begin{align*}
&-\overset{d}{\underset{p,q=1}{\sum}} \left[\dfrac{1}{\abs{Y}}\overset{d}{\underset{k=1}{\sum}}   \displaystyle \int_{Y_{i}}  \left(\widetilde{\mathbf{m}}_i^{pk} \dfrac{\pt \chi_{i}^q}{\pt y_k}(y)+\widetilde{\mathbf{m}}_i^{pq}\right) dy \right]\dfrac{\pt^2 u_{i,0}}{\pt x_p \pt x_q}
\\& \qquad\qquad =- \dfrac{\abs{\Gamma^{y}}}{\abs{Y}} \left( \pt_{t} v_0+\I_{ion}(v_0,{w}_0)-\I_{app}\right)
\end{align*} 
\end{enumerate}
 Finally, we deduce the \textbf{homogenized} equation satisfied by $u_{i,0}$ for the intracellular problem:
\begin{equation}
\B_{xx} u_{i,0}= -\mu_{m} \left(\pt_{t} v_0+\I_{ion}(v_0,{w}_0)-\I_{app}\right) \ \text{on} \ \Omega_T,  
 \label{pbhomin}
 \end{equation}
 where $\mu_{m}=\abs{\Gamma^{y}}/\abs{Y}$. Here, the homogenized operator $\B_{xx}$ (with respect to $y$ and $z$) is defined by:
 \begin{equation*}
 \B_{xx}=-\nabla_x \cdot \left( \doublewidetilde{\mathbf{M}}_i \nabla_x \right)=-\overset{d}{\underset{p,q=1}{\sum}}\dfrac{\pt}{\pt x_p}\left(\doublewidetilde{\mathbf{m}}^{pq}_{i}\dfrac{\pt}{\pt x_q}\right)
 \label{Bxxhomi}
 \end{equation*}
 with the coefficients of the homogenized conductivity matrix $\doublewidetilde{\mathbf{M}}_i=\left( \doublewidetilde{\mathbf{m}}^{pq}_i\right)_{1\leq p,q \leq d}$ defined by:
\begin{equation}
\begin{aligned}
\doublewidetilde{\mathbf{m}}^{pq}_i & :=\dfrac{1}{\abs{Y}}\overset{d}{\underset{k=1}{\sum}}\displaystyle \int_{Y_{i}}  \left(\widetilde{\mathbf{m}}_i^{pk} \dfrac{\pt \chi_{i}^q}{\pt y_k}(y)+\widetilde{\mathbf{m}}_i^{pq}\right) dy
\\ & = \dfrac{1}{\abs{Y}} \dfrac{1}{\abs{Z}}\overset{d}{\underset{k,\ell=1}{\sum}}\displaystyle \int_{Y_{i}} \int_{Z_{c}}    \left[\displaystyle \left(\mathrm{m}^{pk}_{i}+ \mathrm{m}^{p\ell}_{i}\dfrac{\pt \theta_{i}^{k}}{\pt z_\ell}\right) \dfrac{\pt \chi_{i}^q}{\pt y_k}(y)+\left(\mathrm{m}^{pq}_{i}+ \mathrm{m}^{p\ell}_{i}\dfrac{\pt \theta_{i}^{q}}{\pt z_\ell}\right)\right] \ dzdy
\end{aligned}
\label{Mb_i}
\end{equation} 
with the coefficients of the conductivity matrix $\widetilde{\mathbf{M}}_i=\left(\widetilde{\mathbf{m}}^{pk}_i\right)_{1\leq p,k \leq d}$ defined by \eqref{Mtt_i}.
 
\begin{rem} The authors in \cite{henri} treated the initial problem with the coefficients $\mathrm{m}^{pq}_{j}$ depending only on the variable $y$ for $j=i,e$. Using the same two-scale technique, we found three systems to solve and then obtained its homogenized model with respect to $y$ which is well defined in Section \ref{asymextra}. But in the intracellular problem, the coefficients $\mathrm{m}^{pq}_{i}$ depend on two variables $y$ and $z.$ Using a new three-scale expansion method, we obtain nine systems to solve in order to find the homogenized model \eqref{pbhomin} of the initial problem \eqref{pbinintra}. Obtaining this homogenized problem is described in six steps. First, the first five steps help to find the explicit forms of the associated solutions. Second, the last step describes the two-level homogenization whose the coefficients $\doublewidetilde{\mathbf{m}}^{pq}_i$ of the homogenized conductivity matrix $\doublewidetilde{\mathbf{M}}_i$  are integrated with respect to $z$ and then with respect to $y.$ Finally, we obtain the homogenized model defined on $\Omega.$
  \end{rem}
\subsection{Macroscopic Bidomain Model}\label{macro_bid_asym} At macroscopic level, the heart domain coincides with the intracellular and extracellular ones, which  are inter-penetrating and superimposed connected at each point by the cardiac cellular membrane. The homogenized model of the microscopic bidomain model are recuperated from the extracellular and intracellular homogenized equations \eqref{pbhomex}-\eqref{pbhomin}, which is called the macroscopic bidomain model (Reaction-Diffusion system):
\begin{equation}
\begin{aligned}
\mu_{m}\pt_{t} v+\nabla \cdot\left(\widetilde{\mathbf{M}}_{e}\nabla u_{e}\right) +\mu_{m}\I_{ion}(v,w) &= \mu_{m}\I_{app} &\text{ in } \Omega_{T},
\\ \mu_{m}\pt_{t} v-\nabla \cdot\left( \doublewidetilde{\mathbf{M}}_{i}\nabla u_{i}\right)+\mu_{m}\I_{ion}(v,w) &= \mu_{m} \I_{app} &\text{ in } \Omega_{T},
\\ \pt_{t} w-H(v,w) &=0 & \text{ on }  \Omega_{T},
\label{pb_macro}
\end{aligned}
\end{equation}
completed with no-flux boundary conditions on $u_i, u_e$ on $\pt_{\text{ext}} \Omega:$ 
\begin{equation*}
\left(\widetilde{\mathbf{M}}_{e}\nabla u_{e}\right)\cdot\mathbf{n}=\left( \doublewidetilde{\mathbf{M}}_{i}\nabla u_{i}\right) \cdot\mathbf{n}=0 \ \text{ on } \  \Sigma_T:=(0,T)\times\pt_{\text{ext}} \Omega,
\end{equation*}
where $\mathbf{n}$ the outward unit normal to the boundary of $\Omega,$ and by assigning the initial Cauchy condition for the transmembrane potential $v$ and the gating variable $w :$ 
\begin{equation}
v(0,x)=v_0(x)\qquad\text{and}\qquad w(0,x)=w_0(x), \text{ a.e. on } \Omega.
\label{cond_ini}
\end{equation}
 Herein, the conductivity matrices $\widetilde{\mathbf{M}}_{e}$ and $\doublewidetilde{\mathbf{M}}_{i}$ are defined respectively in \eqref{Mb_e}-\eqref{Mb_i}. System  \eqref{pb_macro}-\eqref{cond_ini} corresponds to the sought macroscopic equations. Finally, note that we close the problem by the normalization condition on the extracellular potential for almost all $t\in [0,T],$
\begin{equation*}
\int_{\Omega} u_{e} (t,x)dx=0.
\end{equation*}

\begin{rem} Following \cite{ben,doina}, it is easy to verify that these homogenized conductivity tensors are symmetric, positive definite. Moreover, the functions $\I_{ion}$ and $H(v,w)$ preserve the same form of Fitzhugh-Nagumo model defined in \eqref{ionic_model_asym}.
\end{rem}

\section{Conclusion}
 Many biological and physical phenomena arise in highly heterogeneous media, the properties of which vary on three (or more) length scales. 
In this paper, a new three-scale asymptotic homogenization technique have been established for predicting the bioelectrical behaviors of the cardiac tissue with multiple small-scale configurations. Furthermore, we have presented the main mathematical models to describe the bioelectrical activity of the heart, from the microscopic activity of ion channels of the cellular membrane to the macroscopic properties in the whole heart. We have described how reaction-diffusion systems can be derived from microscopic models of cellular aggregates by homogenization method and a new three-scale asymptotic expansion. 
  
  The present study has some limitations and is open to several improvements. For example, analytical formulas have been found for an ideal particular geometry at the mesoscale and microscale. Nevertheless, the natural next step is to consider more realistic geometries by solving the appropriate cellular problems analytically and numerically.
 
 A key assumption underlying the whole method is periodicity of the micro-structure at both structural levels. This assumption can be considered realistic for specific types of micro-structures only. However, our framework is extended to more complex geometries by taking account of two parameters scaling dependent on the cell geometry on the macroscale. A special attention to the boundary conditions for the unit cell to ensure periodicity. 
 
 The homogenization process described in this work is also suitable for regions far enough from the boundary so that its effect is not felt (for example composite material). To account properly the homogenization process on bounded domains, so-called the boundary-layer technique established by A. Benssousan et al. \cite{ben} could be used (see also the work of Panasenko \cite{pana}). We known results from reiterated and two-scale asymptotic homogenization techniques as particular cases of the proposed method. \\ \\
\textbf{Funding.} This research was supported by IEA-CNRS in the context of HIPHOP project.

\bibliographystyle{plain}
 \bibliography{Hom}

\appendix

\section{Periodic Sobolev space}\label{appA}
In this section, we give the properties which play an important role in the theory of homogenization (see \cite{doina}). For more details on functional analysis, the reader is referred to the following references:  \cite{rudin}, \cite{adams}, \cite{edwards}, \cite{brezis}, \cite{zeidler}.
We denote by  $\mathcal{O}$ the interval in  $\R^d$ defined by :
\begin{equation}
 \mathcal{O}=\left]0,\ell_{1}\right[ \times \cdots \times \left]0,\ell_{d}\right[,
 \label{O}
 \end{equation}
where $\ell_{1}, \dots ,\ell_{d}$ are given positive numbers. We will refer to $\mathcal{O}$ as the reference cell. \\

 We define now the periodicity for functions which are defined almost everywhere.
\begin{defi}Let $\mathcal{O}$ the reference cell defined by \eqref{O} and $f$ a function defined a.e on $\R^d.$\\ The function $f$ is called $\textbf{y}$-periodic, if and only if,
$$f(y+k\ell_{i}e_{i})=f(y) \ \text{p.p. on} \ \R^d,\ \forall k\in \mathbb{Z}, \ \forall i \in \lbrace 1, \dots, d \rbrace, $$
where $\lbrace e_{1}, \dots ,e_{d} \rbrace$ is the canonical basis of $\R^d.$
\label{Y_per}
\end{defi}

\begin{defi}Let $\alpha, \beta \in \R,$ such that $0<\alpha<\beta.$ We denote by $M(\alpha,\beta,\mathcal{O})$ the set of the $d \times d$ matrices $\mathrm{M}=(\mathrm{m}^{pq})_{1 \leq p,q\leq d}\in L^\infty(\mathcal{O})^{d\times d}$ such that :
\begin{equation}
\begin{cases}
(\mathrm{M}(x)\lambda,\lambda)\geq \alpha\abs{\lambda}^2, \\ \abs{\mathrm{M}(x)\lambda}\leq \beta \abs{\lambda},
 \end{cases}
 \end{equation}
for any $\lambda\in \R^d$ and almost everywhere on $ \mathcal{O}.$
\label{M}
\end{defi}

In this part, we introduce a notion of periodicity for functions in the Sobolev space $H^{1}.$ In the sequel, we take $\mathcal{O}$ an open bounded set in $\R^d.$

\begin{defi} Let $C^{\infty}_{per}(\mathcal{O})$ be the subset of $C^{\infty}(\R^d)$ of periodic functions. We denote by $H^{1}_{per}(\mathcal{O})$ the closure of $C^{\infty}_{per}(\mathcal{O})$ for the $H^{1}$-norm, namely, $$H^{1}_{per}(\mathcal{O})=\overline{C^{\infty}_{per}(\mathcal{O})}^{H^{1}(\mathcal{O})}.$$
\end{defi}

\begin{prop}Let $u \in H^{1}_{per}(\mathcal{O}).$ Then $u$ has the same trace on the opposite faces of $\mathcal{O}.$
\label{cond_Y_per}
\end{prop}

In the sequel, we will define the quotient space  $H^{1}_{per}(\mathcal{O})/\R$ and introduce some properties on this space.

\begin{defi}The quotient space $\W_{per}(\mathcal{O})$ is defined by:
$$\W_{per}(\mathcal{O})=H^{1}_{per}(\mathcal{O})/ \R .$$
It is defined as the space of equivalence classes with respect to the following relation:
$$u \simeq v \Leftrightarrow u-v \ \text{is a constant,} \ \forall u,v \in H^{1}_{per}(\mathcal{O}).$$
We denote by $\dot{u}$ the equivalence class represented by $u.$  
\label{w_per}
\end{defi}
\begin{prop}The following quantity:
$$\norm {\dot{u}}_{\W_{per}(\mathcal{O})}=\norm {\nabla u}_{L^2(\mathcal{O})}, \forall u \in \dot{u}, \dot{u} \in \W_{per}(\mathcal{O}) $$ defines a norm on $\W_{per}(\mathcal{O}).$ \\
Moreover, the dual space $(\W_{per}(\mathcal{O}))'$ can be identified with the set:
$$(\W_{per}(\mathcal{O}))'=\lbrace F \in (H^{1}_{per}(\mathcal{O}))' \ \text{tel que} \ F(c)=0,  \ \forall c\in \R \rbrace,$$
with $$F(u)=(F,\dot{u})_{(\W_{per}(\mathcal{O}))',\W_{per}(\mathcal{O})}=(F,u)_{(H^{1}_{per}(\mathcal{O}))',H^{1}_{per}(\mathcal{O})},\ \forall u \in \dot{u}, \dot{u} \in \W_{per}(\mathcal{O}).$$
\end{prop}

\begin{rem}In particular, we can choose a representative element  $u$ of the equivalence class $\dot{u}$ by fixing the constant. Then, we define a particular space of periodic functions with a null mean value as follows:
\begin{equation}
W_{per}(\mathcal{O})=\lbrace u \in H^1_{per}(\mathcal{O}) \ \text{such that} \ \M_{\mathcal{O}}(u)=0\rbrace.
\label{W}
\end{equation}
with $\M_{\mathcal{O}}(u)=\dfrac{1}{\abs{\mathcal{O}}}\displaystyle\int_{\mathcal{O}} u \ dx.$
Its dual coincides with the dual space $ (\W_{per}(\mathcal{O}))'$ and the duality bracket is defined by: $$F(v)=(F,v)_{(W_{\#}(\mathcal{O}))',W_{\#}(\mathcal{O})}=(F,u)_{(H^{1}_{per}(\mathcal{O}))',H^{1}_{per}(\mathcal{O})},\ \forall u \in W_{\#}(\mathcal{O}).$$
Furthermore, by the Poincaré-Wirtinger's inequality, the Banach space $W_{\#}(\mathcal{O})$ has the following norm:
$$\norm {u}_{W_{\#}(\mathcal{O})}=\norm {\nabla u}_{L^2(\mathcal{O})}, \forall u \in W_{\#}(\mathcal{O}).$$
\label{rem_w}
\end{rem}

 In the sequel, we will introduce some elliptic partial differential equations  with different boundary conditions: Neumann and periodic conditions. In these cases, to prove existence and uniqueness, the Lax-Milgram theorem will be applied. Few works are available in the literature about boundary value problems, we cite for instance \cite{lions}, \cite{lions1969}. In this part, we will treat the following partial equation: 
\begin{equation*}
\A u=f \ \text{in} \ \mathcal{O},
\end{equation*} 
with the operator $\A$ defined by: 
\begin{equation}
\A=-\nabla \cdot(\mathrm{M} \nabla)
\label{A}
\end{equation}
where the matrix $\mathrm{M}=(\mathrm{m}^{pq})_{1 \leq p,q \leq d} \in M(\alpha,\beta,\mathcal{O})$ is given by Definition \ref{M} but with different boundary conditions:
\begin{enumerate}
\item[$\bullet$]$\textbf{Nonhomogenous Neumann condition:}$ 
\begin{equation*}
\mathrm{M} \nabla u\cdot n =g \ \text{on} \ \pt \mathcal{O}.
\end{equation*}
\item[$\bullet$]$\textbf{Periodic-Neumann condition:}$ Let $\mathcal{O}_j$ a portion of a reference cell $\mathcal{O}$  given by \eqref{O}, with a boundary $\Gamma$ separate the two regions $\mathcal{O}_j$ and $\mathcal{O} \setminus \mathcal{O}_j$. So, we have : $$\pt \mathcal{O}_j=(\pt \mathcal{O} \cap \pt \mathcal{O}_j)\cup \Gamma.$$ The boundary condition which plays an essential role in the homogenization of perforated periodic media, namely,
\begin{equation*}
\begin{cases}
u  \ \ \textbf{y}\text{-periodic}, \\ \mathrm{M} \nabla u\cdot n =g \ \text{on} \ \Gamma.
\end{cases}
\end{equation*}
\end{enumerate}

\begin{thm}$(\text{Nonhomogenous Neumann condition})$\\
We consider the following problem:
\begin{equation}
 \begin{cases}
\A u=f \ \text{in} \ \mathcal{O},\\ \mathrm{M} \nabla u\cdot n =g \ \text{on} \ \pt \mathcal{O}. 
 \end{cases}
 \label{AuNeu}
 \end{equation}
  with the operator $\A$ defined by \eqref{A}. Its variational formulation is:
 \begin{equation}
 \begin{cases}
\text{Find} \ u \in H^{1}(\mathcal{O}) \ \text{such that} \\ a_{\mathcal{O}}(u,v)=(f,v)_{H^{-1}(\mathcal{O}),H^{1}(\mathcal{O})}+(g,v)_{H^{-\frac{1}{2}}(\pt\mathcal{O}),H^{\frac{1}{2}}(\pt\mathcal{O})} \ \forall v \in H^{1}(\mathcal{O}),
 \end{cases}
 \label{FvAuNeu}
\end{equation}
with $a_{\mathcal{O}}$ defined by: 
\begin{equation*}
a_{\mathcal{O}}(u,v)=\int_{\mathcal{O}} \mathrm{M} \nabla u \nabla v dx, \ \forall u,v \in H^{1}(\mathcal{O}).
\end{equation*} 
  We take $V=H^{1}(\mathcal{O}).$ Suppose that $f\in L^2(\mathcal{O})$ and $g \in H^{\frac{1}{2}}(\pt \mathcal{O})$ satisfy the following compatibility condition:
  \begin{equation}
  (f,1)_{H^{-1}(\mathcal{O}),H^{1}(\mathcal{O})}+(g,1)_{H^{-\frac{1}{2}}(\pt\mathcal{O}),H^{\frac{1}{2}}(\pt\mathcal{O})}=0.
\end{equation}    

 Then, the problem \eqref{AuNeu}-\eqref{FvAuNeu} has a unique solution $u \in H^{1}(\mathcal{O})$. Moreover,
\begin{equation*}
{\norm {u}_{H^{1}(\mathcal{O})}}\leq \dfrac{1}{\alpha_0}\left(\norm {f}_{L^2(\mathcal{O})}+C_{\gamma}\norm {g}_{H^{-\frac{1}{2}}(\pt \mathcal{O})}\right) ,
\end{equation*}
where $\alpha_0=\min(1,\alpha)$ and $C_\gamma$ is the trace constant.
\label{cond_neu}
\end{thm}

\begin{thm}$(\text{Periodic-Newmann condition})$\\
Let $\mathcal{O}_j$ a portion of a unit cell $\mathcal{O}$  given by \eqref{O}, with Lipschitz continuous boundary $\Gamma$ separate the two regions $\mathcal{O}_j$ and $\mathcal{O} \setminus \mathcal{O}_j$. Consider the following problem:
\begin{equation}
 \begin{cases}
\A u=f \ \text{in} \ \mathcal{O}_j,\\ u  \ y\text{-periodic}, \\ \mathrm{M} \nabla u\cdot n =g \ \text{on} \ \Gamma.  
 \end{cases}
 \label{Auperneu}
 \end{equation}
We take $V=\W_{per}(\mathcal{O}_j).$ Then, for any $f\in(\W_{per}(\mathcal{O}_j))'$ and for any $g \in H^{\frac{1}{2}}(\Gamma),$ the variational formulation of the problem \eqref{Auperneu} is:
 \begin{equation}
 \begin{cases}
\text{Find} \ \dot{u} \in \W_{per}(\mathcal{O}_j) \ \text{such that} \\ \dot{a}_{\mathcal{O}_j}(\dot{u},\dot{v})=(F,\dot{v})_{(\W_{per}(\mathcal{O}_j))',\W_{per}(\mathcal{O}_j)} \ \forall \dot{v}\in \W_{per}(\mathcal{O}_j),
 \end{cases}
 \label{FvAuperneu}
\end{equation}
 
with $a_{\mathcal{O}_j}$ is given by: 
\begin{equation*}
\dot{a}_{\mathcal{O}_j}(\dot{u},\dot{v})=\int_{\mathcal{O}_j} \mathrm{M} \nabla u \nabla v dy, \ \forall u \in \dot{u},\ \forall v \in \dot{v},
\end{equation*}
and $F$ is defined by:
\begin{equation*}
(F,\dot{v})_{(\W_{per}(\mathcal{O}_j))',\W_{per}(\mathcal{O}_j)}=\displaystyle \int_{\Gamma}\mathrm{M}_i \nabla u \cdot n v \ d\sigma_y+\displaystyle \int_{\mathcal{O}_{j}} fv \ dy, \ \forall v \in\dot{v}, \ \forall \dot{v} \in \W_{per}(\mathcal{O}_j),
\end{equation*}
where $n$ denotes the unit outward normal to $\Gamma$.
 
 Assume that $\mathrm{M}$ belongs to $M(\alpha,\beta,\mathcal{O})$ with $y$-periodic coefficients. Suppose that $F$ belongs to $(\W_{per}(\mathcal{O}_j))'$ which equivalent to
 \begin{equation*}
 (F,1)_{(\W_{per}(\mathcal{O}_j))',\W_{per}(\mathcal{O}_j)}=0.
\end{equation*}  
Then problem \eqref{FvAuperneu} has a unique weak solution. Moreover, we have the following estimation:
 \begin{equation*}
{\norm {\dot{u}}_{\W_{per}(\mathcal{O}_j)}}\leq \dfrac{1}{\alpha_0}\left(\norm {f}_{L^2(\mathcal{O}_j)}+C_{\gamma}\norm {g}_{H^{-\frac{1}{2}}(\Gamma)}\right).
\end{equation*}
\label{cond_per_neu}
where $\alpha_0=\min(1,\alpha)$ and $C_\gamma$ is the trace constant.
\end{thm} 

  By the definition of $\W_{per},$ the previous theorem shows that the problem \eqref{Auperneu} admits a solution in $H^1_{per},$ defined up to an additive constant. If we take the particular case $V=W_{\#}(\mathcal{O})$ defined by \eqref{W}, we obtain the same result.

\end{document}